\documentclass{amsart}

\usepackage{amssymb}
\usepackage{color}
\usepackage{graphicx}

\newcommand{\dx}{\partial_x}
\newcommand{\dz}{\partial_z}
\newcommand{\dt}{\partial_t}
\newcommand{\eps}{\varepsilon}
\newcommand{\R}{\mathbb R}
\newcommand{\N}{\mathbb N}
\newcommand{\cst}{\mbox{\textnormal{Cst }}}
\newtheorem{prop}{Proposition}
\newtheorem{theo}{Theorem}
\newtheorem{defi}{Definition}
\newtheorem{rema}{Remark}

\title[The hydrodynamical relevance of the CH and DP equations]{The hydrodynamical relevance of the Camassa-Holm and Degasperis-Procesi equations}
\author{Adrian Constantin and David Lannes}
\address{School of Mathematics, Trinity College, Dublin 2, Ireland}
\email{adrian@maths.tcd.ie}
\address{Universit\'e Bordeaux I; IMB and CNRS UMR 5251, 351 Cours de la Lib\'eration, 33405 Talence Cedex, France}
\email{lannes@math.u-bordeaux1.fr}

\begin{document}

\maketitle

\begin{abstract}
In recent years two nonlinear dispersive partial differential equations have attracted a lot of attention 
due to their integrable structure. We prove that both equations arise in the modeling of 
the propagation of shallow water waves over a flat bed. The equations capture stronger nonlinear effects  than 
the classical nonlinear dispersive Benjamin-Bona-Mahoney and Korteweg-de Vries equations. In particular, they
accomodate wave
breaking phenomena.
\end{abstract}

\section{Introduction}

The study of water waves is a fascinating subject because the phenomena are familiar and the 
mathematical problems are various cf. \cite{Wh}. Due to the relative intractability of the governing 
equations for water waves in regard to inferring from their direct study 
qualitative or quantitative conclusions about the propagation of waves at the water's surface, 
from the earliest days in the development of hydrodynamics many competing models were 
suggested. Until the second half of the 20'th century, the study of water waves was 
confined almost exclusively to linear theory \cite{Cr}. While linearisation gives insight for small 
perturbations on water initially at rest, its applicability fails for waves that are not 
small perturbations of a flat water surface. For example, linear water wave theory 
gives no insight into the study of phenomena which are manifestations of genuine nonlinear 
behaviour, like breaking waves breaking and solitary waves \cite{St}. Many nonlinear models 
for water waves have been suggested to capture the existence of solitary water waves and 
the associated phenomenon of soliton manifestation \cite{Jo-b}. The most prominent example 
is the Korteweg-de Vries  (KdV) equation \cite{KdV}, the only member of the wider family of  
BBM-type equations \cite{BBM} that is integrable and relevant for the phenomenon of soliton manifestation \cite{DJ}.
Another development of models for water waves was initiated in order to gain insight into wave 
breaking, one of the most fundamental aspects of water waves for which there appears 
to be no satisfactory mathematical theory \cite{Wh}. Starting from the observation that the strong 
dispersive effect incorporated into the KdV model prevents wave breaking, Whitham (see 
the discussion in \cite{Wh}) initiated the quest for equations that are simpler than the governing 
equations for water waves and which could model breaking waves. The physical validity of the 
first proposed models is questionable but two recently derived nonlinear integrable equations, 
the Camassa-Holm equation \cite{CH} and the Degasperis-Procesi equation \cite{DP}, 
possess smooth solutions that develop singularities in finite time via a process that 
captures the essential features of breaking waves cf. \cite{Wh}: the solution 
remains bounded but its slope becomes unbounded. Our aim is to prove the 
relevance of these two equations as models for the propagation of shallow water waves, 
proving that both are valid approximations to 
the governing equations for water waves. In our investigation we put earlier (formal) 
asymptotic procedures due to Johnson \cite{Jo} on a firm and mathematically rigorous basis. 
We also investigate in what sense these two models give us insight into the wave breaking 
phenomenon by some simple numerical
computations.

\subsection{Unidirectional asymptotics for water waves}

For one dimensional surfaces, the water waves equations read, in nondimensionalized form,
\begin{equation}\label{ww}
	\left\lbrace
	\begin{array}{lcl}
	\mu\dx^2\Phi+\dz\Phi^2=0 &\mbox{ in }&\Omega_t,\\
	\dz\Phi=0,&\mbox{ at }& z=0,\\
	\dt \zeta-\frac{1}{\mu}(-\mu\dx\zeta\dx\Phi+\dz\Phi)=0
	&\mbox{ at }& z=\eps\zeta,\\
	\dt\Phi+\frac{\eps}{2}(\dx\Phi)^2+\frac{\eps}{2\mu}(\dz\Phi)^2=0&\mbox{ at }& z=\eps\zeta,
	\end{array}\right.
\end{equation}
where $x\mapsto \eps\zeta(t,x)$ parameterizes the elevation of the free surface at time $t$, $\Omega_t=\{(x,z),-1<z<\eps \zeta(t,x)\}$ is the fluid domain
delimited by the free surface and the flat bottom $\{z=-1\}$, and where
$\Phi(t,\cdot)$ (defined on $\Omega_t$) is the velocity potential
associated to the flow (that is, the two-dimensional 
velocity field ${\bf v}$ is given by ${\bf v}=(\dx\Phi,\dz\Phi)^T$).
Finally, $\eps$ and $\mu$ are two dimensionless parameters defined as
$$
	\eps=\frac{a}{h},\qquad \mu=\frac{h^2}{\lambda^2},
$$
where $h$ is the mean depth, $a$ is the typical amplitude and $\lambda$ the typical wavelength of the waves under consideration. Making assumptions
on the respective size of $\eps$ and $\mu$, one is led to derive
(simpler) asymptotic models from (\ref{ww}).\\
In the shallow-water scaling ($\mu\ll1$), one can derive the so-called
Green-Naghdi equations (see \cite{GN} for the derivation, and \cite{AL} for
a rigorous justification), without any assumption on 
$\eps$ (that is, $\eps=O(1)$). For one dimensional surfaces and flat
bottoms, these equations couple the free surface
elevation $\zeta$ to the vertically averaged horizontal component of the
velocity,
\begin{equation}\label{averaged}
	u(t,x)=\frac{1}{1+\eps\zeta}\int_{-1}^{\eps\zeta}\dx \phi(t,x,z)dz;
\end{equation}
and can be written as
\begin{equation}\label{eq1}
	\left\lbrace
	\begin{array}{l}
	\zeta_t+\big[(1+\eps\zeta)u\big]_x=0\\
	u_t+\zeta_x+\eps u u_x =\frac{\mu}{3}\frac{1}{1+\eps\zeta}
	\big[(1+\eps\zeta)^3(u_{xt}+\eps u u_{xx}-\eps u_x^2)\big]_x,
	\end{array}\right.
\end{equation}
where $O(\mu^2)$ terms have been discarded.\\
If we make the \emph{additionnal} assumption that $\eps\ll1$, then
the above system reduces at first order to a wave equation of speed
$\pm 1$ and any perturbation of the surface splits up into two
components moving in opposite directions. A natural issue is therefore
to describe more accurately the motion of these two ``unidirectional''
waves. In the
so-called long-wave regime
\begin{equation}\label{scalingLW}
	\mu\ll 1,\qquad \eps=O(\mu),
\end{equation}
Korteweg and de Vries \cite{KdV} found that, say, the right-going wave
should satisfy the KdV equation
$$
	u_t+u_x+\eps \frac{3}{2}uu_x+\mu\frac{1}{6}u_{xxx}=0
$$
(and $\zeta=u+O(\eps,\mu)$), which at leading order reduces to the expected
transport equation at speed $1$. More recently, it has been
noticed by Benjamin, Bona, Mahoney \cite{BBM} that the KdV equation belongs
to a wider class of equations (the BBM equations, first used by Peregrine \cite{Pe} 
and sometimes also called the regularized long-wave equations) which provide an
approximation of the exact water waves equations of the same accuracy
as the KdV equation:
\begin{equation}\label{BBM}
	u_t+u_x+\frac{3}{2}\eps u u_x+
	\mu (\alpha u_{xxx}+\beta u_{xxt})
	=0,\quad\mbox{ with }\quad \alpha-\beta=\frac{1}{6}.
\end{equation}
The equations (\ref{BBM}) contain both nonlinear
effects (the $uu_x$ term) and dispersive effects (the $u_{xxx}$and $u_{xxt}$ terms) due to the scaling (\ref{scalingLW}). However, these equations
do not account correctly for large amplitude waves, whose behavior is
more nonlinear than dispersive. For such waves, characterized by
larger values of $\eps$, it is natural to investigate the following
scaling (which we call Camassa-Holm scaling):
\begin{equation}\label{scalingCH}
	\mu\ll1,\qquad \eps=O(\sqrt{\mu}).
\end{equation}
With this scaling, one still has $\eps\ll1$ and thus the same reduction
to a simple wave equation at leading order; the dimensionless parameter
is however larger here than in the long wave scaling, and the nonlinear
effects are therefore stronger. In particular, a stronger nonlinearity could allow 
the appearance of breaking waves - a fundamental phenomenon in the 
theory of water waves that is not captured by the BBM equations. We show in this paper that
the correct generalization of the BBM equations (\ref{BBM}) under
the scaling (\ref{scalingCH}) is provided by  
the following class of equations:
\begin{equation}\label{CHfamily}
	u_t+u_x+\frac{3}{2}\eps u u_x+
	\mu (\alpha u_{xxx}+\beta u_{xxt})
	=\eps\mu(\gamma u u_{xxx}+\delta u_x u_{xx})
\end{equation}
(with some conditions on $\alpha$, $\beta$, $\gamma$ and $\delta$).\\
Notice that for an 
equation of the family (\ref{CHfamily}) to be well-posed it is necessary that $\beta \le 0$, as 
one can see by analyzing the linear part via Fourier transforms. 
 
We want to insist on the fact that (\ref{CHfamily}) provides an 
approximation of the same order $O(\mu^2)$ as the BBM equations
(\ref{BBM}) to the Green-Naghdi equations. The only difference in the
derivation of these equations lies in the different scalings
(\ref{scalingLW}) and (\ref{scalingCH}). Since of course $O(\mu)=O(\sqrt{\mu})$
when $\mu$ is small, the long-wave scaling (\ref{scalingLW}) is contained
in the CH scaling (\ref{scalingCH}), and consequently, the BBM equations
can be recovered as a specialization of (\ref{CHfamily}) when 
$\eps=O(\mu)$ and not only $O(\sqrt{\mu})$.

\subsection{The Camassa-Holm and Degasperis-Procesi equations}

Among the various type of equations (\ref{CHfamily}) with $\beta \le 0$ there are only 
two with a bi-Hamiltonian structure:  the Camassa-Holm and the Degasperis-Procesi equations 
\cite{Iv}. Notice that while the KdV equation has a bi-Hamiltonian structure (see \cite{DJ}), this 
is not the case for the other members of the BBM family of equations (\ref{BBM}). The importance of 
a bi-Hamiltonian structure lies in the fact that in general it represents the hallmark of a completely 
integrable Hamiltonian system whose solitary wave solutions are solitons, that is, localized 
waves that recover their shape and speed after interacting nonlinearly with another wave of the 
same type (see \cite{DJ, Jo-b}).

\subsubsection{Camassa-Holm equations}

The Camassa-Holm (CH) equations are usually written under the form
\begin{equation}\label{CHeq}
	U_t+\widehat{\kappa}U_x+3UU_x-U_{txx}=2U_xU_{xx}+UU_{xxx},
\end{equation}
with $\widehat{\kappa}\in \R$. A straightforward scaling argument shows
that if $\widehat{\kappa}\neq 0$, 
(\ref{CHeq}) can be written under the form (\ref{CHfamily})
by setting $u(t,x)=aU(b(x-vt),ct)$ and $a=\frac{2}{\eps\widehat{\kappa}}$,
$b^2=-\frac{1}{\beta\mu}$, $v=\frac{\alpha}{\beta}$, $c=\frac{b}{\widehat{\kappa}}(1-v)$
(which requires $\beta<0$ and leads to $\gamma=-\frac{\beta}{2}$
and $\delta=2\gamma$). This motivates the following definition:
\begin{defi}\label{defiCH}
	We say that (\ref{CHfamily}) is a Camassa-Holm equation if
	the following conditions hold:
	$$
	\beta<0,\quad \alpha\neq \beta, \quad\beta=-2\gamma,\quad 
	\delta=2\gamma.
	$$
	 For all $\widehat{\kappa}\neq 0$, 
	the solution $u$ to (\ref{CHfamily}) is transformed
	into a solution $U$ to (\ref{CHeq}) by the transformation
	$$
	U(t,x)=\frac{1}{a}u(\frac{x}{b}+\frac{v}{c}t,\frac{t}{c}),
	$$
	with
	$a=\frac{2}{\eps\widehat{\kappa}}(1-v)$,
	$b^2=-\frac{1}{\beta\mu}$,
	$v=\frac{\alpha}{\beta}$,
	and $c=\frac{b}{\widehat{\kappa}}(1-v)$.
\end{defi}

First derived as a bi-Hamiltonian system by Fokas \& Fuchssteiner \cite{FF}, the 
equation (\ref{CHeq}) gained prominence after Camassa-Holm \cite{CH} independently 
re-derived it as an approximation to the Euler equations of hydrodynamics and 
discovered a number of the intriguing properties of this equation. In \cite{CH} a  
Lax pair formulation of (\ref{CHeq}) was found, a fact which lies at the core of showing 
via direct and inverse scattering \cite{Co, CGI} that (\ref{CHeq}) is a completely integrable 
Hamiltonian system: for a large class of initial data, solving (\ref{CHeq}) amounts 
to integrating an infinite number of linear first-order ordinary differential equations which 
describe the evolution in time of the action-angle variables. The Camassa-Holm equations shares 
with KdV this integrability property as well as the fact that its solitary waves are solitons \cite{CS, CGI}. 
We refer to \cite{DJ} for a discussion of these properties in the context of the KdV model. 

\subsubsection{Degasperis-Procesi equations}

The Degasperis-Procesi (DP) equations are usually written  under the form
\begin{equation}\label{DPeq}
	U_t+\widehat{\kappa}U_x+4UU_x-U_{txx}=3U_xU_{xx}+UU_{xxx},
\end{equation}
with $\widehat{\kappa}\in \R$. The same scaling arguments as
for the CH equation motivate the following definition:
\begin{defi}\label{defiDP}
	We say that (\ref{CHfamily}) is a Degasperis-Procesi equation if
	the following conditions hold:
	$$
	\beta<0, \quad \alpha\neq \beta,\quad 
	\beta=-\frac{8}{3}\gamma\quad 
	\delta=3\gamma.
	$$ 
	For all $\widehat{\kappa}\neq 0$, 
	the solution $u$ to (\ref{CHfamily}) is transformed
	into a solution $U$ to (\ref{DPeq}) by the transformation
	$$
	U(t,x)=\frac{1}{a}u(\frac{x}{b}+\frac{v}{c}t,\frac{t}{c}),
	$$
	with
	$a=\frac{8}{3\eps\widehat{\kappa}}(1-v)$,
	$b^2=-\frac{1}{\beta\mu}$,
	$v=\frac{\alpha}{\beta}$ and 
	$c=\frac{b}{\widehat{\kappa}}(1-v)$.
\end{defi}

Equation (\ref{DPeq}), first derived in \cite{DP}, is also known to have a Lax pair formulation \cite{DHH} 
and its solitary waves interact like solitons \cite{Ma}. Just like the KdV equation (see \cite{DJ}) and 
the Camassa-Holm equation (see \cite{Le}), the Degasperis-Procesi equation has infinitely many 
integrals of motion.

\subsection{Wave breaking}

In addition to the properties of (\ref{CHeq}) and (\ref{DPeq}) mentioned before, the 
importance of these two equations is enhanced by their relevance to the modeling of 
wave breaking, one of the most important but mathematically still quite elusive phenomena 
encountered in the study of water waves.
 
\begin{defi}\label{defibreaking}
	We say that there is wave breaking for an equation of 
	the form (\ref{CHfamily}),
	if there exists a time $0<t^{\eps,\mu}<\infty$ and 
	solutions $u$ to (\ref{CHfamily}) such that
	$$
	u\in L^\infty([0,t^{\eps,\mu}]\times\R)
	\quad\mbox{ and }\quad
	\lim_{t\to t^{\eps,\mu}}\vert \dx u(t,\cdot)\vert_{\infty}=\infty.
	$$
\end{defi}

The rationale of this definition is that if the flow velocity $u$ is to first order a good approximation 
to the surface profile $\zeta$, it is then reasonable to expect that the above blow-up pattern has a 
similar counterpart in terms of $\zeta$. If this were the case, the boundedness of the wave height 
in combination with an unbounded slope captures the main features of a breaking wave \cite{Wh}. 

For KdV in particular, as well as for any other member of the BBM family (\ref{BBM}), all smooth 
initial data $u(0,\cdot)$ decaying at infinity develop into solutions defined for all times (see e.g. 
\cite{Ta}) so that the BBM family \cite{BBM} does not model wave breaking \cite{ABLS, SS}. To remedy this 
shortcoming of the KdV equation Whitham proposed to formally replace the dispersive 
term $u_{xxx}$ by a convolution with a singular function chosen so that the newly 
obtained equation presents wave breaking (see \cite{FW, Wh}). 
However, this formal process destroys the integrability and soliton features of the KdV equation. In 
contrast to this, both (\ref{CHeq}) and (\ref{DPeq}) admit breaking waves in the sense of 
Definition \ref{defibreaking} cf. \cite{BC, CE, Mc, Mo} for (\ref{CHeq}) and \cite{ELY} for (\ref{DPeq}). 
In this paper we explore the wave breaking phenomenon for both equations not in the restricted 
sense provided by Definition \ref{defibreaking} but by studying the nonlinear equation describing the evolution in time of the free surface. While our results vindicate the fact that it is appropriate to use in this 
context Definition \ref{defibreaking} to describe breaking waves, there is a slight twist. To be more precise, 
let us distinguish between two types of breaking waves that can be observed. In a plunging breaker 
the slope of the wave approaches $-\infty$ at the breaking location as we reach breaking time, while 
in a surging breaker the slope becomes $+\infty$. Considering the case of the Camassa-Holm equation
(\ref{CHeq}), due to the fact that singularities 
in a smooth solution can appear 
only if $\inf_{x \in \R} \,\{u_x(t,x)\} \to - \infty$ as we approach breaking 
time while $\sup_{x \in \R}\,\{|u(t,x)|\}$ remains uniformly bounded (see \cite{Co}), we would expect 
to observe a plunging breaker at the free surface. However, 
as we recall here, the Camassa-Holm equation describes the
behavior of the vertically averaged horizontal component of the velocity $u$;
the free surface elevation $\zeta$ can be given in terms of $u$
($\zeta=u+\frac{\eps}{4}u^2+\mu\frac{1}{6}u_{xt}+O(\eps\mu)$, see
Proposition \ref{prop1} below) but such an asymptotic expression of course
breaks down when $u_x$ becomes singular, and cannot be used to describe the
behavior of the free surface when there is ``wave breaking'' for the velocity.
Since ``wave breaking'' is a very intuitive notion when it refers to
the free surface elevation, we show in this 
article that it is possible to revert the usual approach; that is, we
derive an evolution equation for the surface elevation $\zeta$ and give
an asymptotic expression for the velocity $u$ in terms of $\zeta$
(cf. \S \ref{sectsurf}). This allows us to prove that wave breaking indeed
occurs for the surface elevation but that, as opposed to what happens for
the velocity, this is a surging breaker! The difference
between plunging and surging
breakers is graphically illustrated by numerical computations in \S \ref{numerics}.

\section{Derivation of asymptotical equations for the unidirectional limit
of the Green-Naghdi equations}

We derive here asymptotical equations to the Green-Naghdi equations in 
the Camassa-Holm scaling (\ref{scalingCH}).
We recall that the Green-Naghdi equations are given by 
$$
	\left\lbrace
	\begin{array}{l}
	\zeta_t+\big[(1+\eps\zeta)u\big]_x=0\\
	u_t+\zeta_x+\eps u u_x =\frac{\mu}{3}\frac{1}{1+\eps\zeta}
	\big[(1+\eps\zeta)^3(u_{xt}+\eps u u_{xx}-\eps u_x^2)\big]_x.
	\end{array}\right.
$$
Since we work under the Camassa-Holm scaling, we 
restrict our attention to values of $\eps$ and $\mu$
satisfying
\begin{equation}\label{cond1}
	(\eps,\mu)\in {\mathcal P}:=\{\mu\in (0,\mu_0),\eps\leq M\sqrt{\mu}\},
\end{equation}
for some $\mu_0>0$ and $M>0$. Equations for the velocity $u$ (including the
CH and DP equations) are first derived in
\S \ref{sectvel}, and equations for the surface elevation $\zeta$ are
obtained in \S \ref{sectsurf}. The considerations we make on the derivation 
of these equations are related to the approach initiated by Johnson \cite{Jo}, 
approach that is substantiated and extended by our analysis. In addition, we 
explore the wave breaking phenomenon.

\subsection{Equations on the velocity}\label{sectvel}

At leading order, the Green-Naghdi equations degenerate into a simple
wave equation of speeds $\pm 1$; including the $O(\eps)$ terms, one can
easily check that the (say) right-going component of the wave must
satisfy
\begin{equation}\label{eq2}
	u_t+u_x+\frac{3}{2}\eps u u_x=0,
\end{equation}
and $\zeta=u+O(\eps)$.\\
If we want to find an asymptotic at order $O(\mu^{2})$ (recall
that $\eps=O(\mu^{1/2})$), it is therefore natural to look for $u$
as a solution of a perturbation of (\ref{eq2}) including
terms of order $O(\mu)$ and $O(\eps\mu)$ similar to those present
in (\ref{eq1}). Thus we want $u$ to solve an equation of the form
(\ref{CHfamily})
where $\alpha$, $\beta$, $\gamma$ and $\delta$ are coefficients to 
be determined. We prove in this section that under certain conditions
on the coefficients, one can associate to the solutions of
(\ref{CHfamily}) a family of approximate solutions
\emph{consistent} with the Green-Naghdi equations (\ref{eq1}) in the
following sense:
\begin{defi}\label{defi1}
	Let $\mu_0>0$, $M>0$, $T>0$ and ${\mathcal P}$ 
	be as defined in (\ref{cond1}).\\
	A family 
	$(\zeta^{\eps,\mu},u^{\eps,\mu})_{(\eps,\mu)\in{\mathcal P}}$
	is \emph{consistent} (of order $s\geq 0$ and on $[0,\frac{T}{\eps}]$) 
	with the Green-Naghdi equations (\ref{eq1}) if
	for all $(\eps,\mu)\in {\mathcal P}$, 
	$$
	\left\lbrace
	\begin{array}{l}
	\zeta_t+\big[(1+\eps\zeta)u\big]_x=\mu^2 r_1^{\eps,\mu}\\
	u_t+\zeta_x+\eps u u_x =\frac{\mu}{3}\frac{1}{1+\eps\zeta}
	\big[(1+\eps\zeta)^3(u_{xt}+\eps u u_{xx}-\eps u_x^2)\big]_x
	+\mu^2 r_2^{\eps,\mu};
	\end{array}\right.
	$$
	with $(r_1^{\eps,\mu},r_2^{\eps,\mu})_{(\eps,\mu)\in{\mathcal P}}$ 
	bounded in 
	$L^\infty([0,\frac{T}{\eps}],H^s(\R)^2)$.
\end{defi}
The following proposition shows that there is a one parameter family
of equations of the form (\ref{CHfamily}) consistent with the Green-Naghdi equations.
\begin{prop} \label{prop1}
	Let $p\in\R$ and asssume that 
	$$
	\alpha=p,\quad \beta=p-\frac{1}{6},\quad
	\gamma=-\frac{3}{2}p-\frac{1}{6},\quad
	\delta=-\frac{9}{2}p-\frac{23}{24}.
	$$
	Then there exists $D>0$ such that:
	\begin{itemize}
	\item For all $s\geq 0$ and $T>0$,
	\item  For all bounded family 
	$(u^{\eps,\mu})_{(\eps,\mu)\in{\mathcal P}}\in 
	C([0,\frac{T}{\eps}];H^{s+D}(\R))$ solving (\ref{CHfamily}),
	\end{itemize}
	the family $(\zeta^{\eps,\mu},u^{\eps,\mu})_{(\eps,\mu)
	\in {\mathcal P}}$,
	with (omitting the indexes $\eps,\mu$)
	$$
	\zeta:=u+\frac{\eps}{4}u^2+\mu\frac{1}{6}u_{xt}
	-\eps\mu\big[\frac{1}{6}uu_{xx}+\frac{5}{48}
	u_x^2\big],
	$$
	is consistent (of order $s$ and on $[0,\frac{T}{\eps}]$) with 
	the Green-Naghdi equations (\ref{eq1}).
\end{prop}
\begin{rema}\label{rm1}
{\bf i.} One can recover the equations (26a) and (26b) of
\cite{Jo} with $p=-\frac{1}{12}$ (and thus $\alpha=-\frac{1}{12}$,
$\beta=-\frac{1}{4}$, $\gamma=-\frac{1}{24}$ and $\delta=-\frac{7}{12}$)
and $p=\frac{1}{6}$ (and thus $\alpha=\frac{1}{6}$, $\beta=0,$ 
$\gamma=-\frac{5}{12}$
and $\delta=-\frac{41}{24}$) respectively.\\
{\bf ii.} The one parameter family of equations (\ref{CHfamily})
considered in the proposition admits only one representant for which
$\delta=2\gamma$ (obtained for
$p=-\frac{5}{12}$, and thus $\gamma=\frac{11}{24},$  $\delta=\frac{11}{12}$). 
However, since $\beta=-\frac{7}{12}\neq -2\gamma$, 
the corresponding equation is not 
a Camassa-Holm equation in the sense of Definition \ref{defiCH}.\\
{\bf iii.} There is no possible choice of $p$ such that $\delta=3\gamma$
in Proposition \ref{prop1}. Consequently, none of this one parameter
family of equations is a Degasperis-Procesi equation. Notice that among all equations  
(\ref{CHfamily}) with $\beta \le 0$ there are only 
two with a bi-Hamiltonian structure:  the Camassa-Holm and the Degasperis-Procesi equations 
\cite{Iv}. 
\end{rema}

\begin{proof} For the sake of simplicity, we use the notation $O(\mu)$, 
$O(\mu^2)$, etc., without explicit mention to the functional normed space
to which we refer. A precise statement has been given in Definition 
\ref{defi1}; it would be straightforward but quite heavy to maintain this
formalism throughout the proof.\\
{\bf Step 1.} If $u$ solves (\ref{CHfamily}) then one also has
\begin{equation}\label{eq4}
	u_t+u_x+\eps\frac{3}{2}uu_x+\mu a u_{xxt}=\eps\mu
	\big[b uu_{xx}+cu_x^2\big]_x+O(\mu^2),
\end{equation}
with $a=\beta-\alpha$, $b=\gamma+\frac{3}{2}\alpha$ and
$c=\frac{1}{2}(\delta+3\alpha-\gamma)$.\\
Differentiating (\ref{CHfamily}) twice with 
respect to $x$, one gets indeed
$$
	u_{xxx}=-u_{xxt}-\frac{3}{2}\eps \dx^2(u u_x)+O(\mu),
$$
and we can replace the $u_{xxx}$ term of (\ref{CHfamily}) by this expression to
get (\ref{eq4}).

\noindent
{\bf Step 2.} We seek $v$ such that if $\zeta=u+\eps v$ and 
$u$ solves (\ref{CHfamily}) then 
the second equation of (\ref{eq1}) is satisfied ut to a $O(\mu^{2})$ term.
This is equivalent to checking that
\begin{eqnarray*}
	u_t+[u+\eps v]_x+\eps u u_x-\frac{\mu}{3}u_{xxt}
	&=&
	\frac{\eps\mu}{3}\big[-uu_{xxt}+[3uu_{xt}+uu_{xx}-u_x^2]_x\big]
	+O(\mu^2),\\
	&=&-\frac{\eps\mu}{3}\big[uu_{xx}+\frac{3}{2}u_x^2\big]_x
	+O(\mu^2)
\end{eqnarray*}
the last line being a consequence of the identity $u_t=-u_x+O(\eps)$ 
provided by (\ref{eq4}). The above equation can be recast
under the form:
\begin{eqnarray*}
	\eps v_x+\big[u_t+u_x+\eps\frac{3}{2}u u_x+\mu a  u_{xxt}
	-\eps\mu \big[b u u_{xx}+c u_x^2\big]_x\big]\\
	=\frac{\eps}{2}u u_x+\mu(a+\frac{1}{3})u_{xxt}
	-\eps\mu \big[(b+\frac{1}{3}) u u_{xx}+(c+\frac{1}{2}) u_x^2\big]_x+O(\mu^2).
\end{eqnarray*}
From Step 1, we know that the term between 
brackets in the lhs of this equation is of order $O(\mu^ 2)$, so that
that the second equation of (\ref{eq1}) is satisfied up to $O(\mu^2)$
terms if
$$
	\eps v_x=\frac{\eps}{2}u u_x+\mu(a+\frac{1}{3})u_{xxt}
	-\eps\mu \big[(b+\frac{1}{3}) u u_{xx}+(c+\frac{1}{2}) u_x^2\big]_x
	+O(\mu^2),
$$
so that we can take
\begin{equation}\label{eq5}
	\eps v=\frac{\eps}{4}u^2+\mu(a+\frac{1}{3})u_{xt}
	-\eps\mu\big[(b+\frac{1}{3})uu_{xx}+(c+\frac{1}{2})
	u_x^2\big].
\end{equation}

\noindent
{\bf Step 3.} We choose the coefficients $\beta$, $\gamma$ and $\mu$ such that
the first equation of (\ref{eq1}) is also satisfied up to $O(\mu^2)$ terms.
This is equivalent to checking that
\begin{equation}\label{eq6}
	[u+\eps v]_t+[(1+\eps u)u]_x+\eps^2[vu]_x=O(\mu^2).
\end{equation}
First remark that one infers from (\ref{eq5}) that
\begin{eqnarray*}
	\eps\dt v&=&\frac{\eps}{2}uu_t+\mu(a+\frac{1}{3})u_{xtt}
	-\eps\mu\big[(b+\frac{1}{3})uu_{xx}+(c+\frac{1}{2})
	u_x^2\big]_t.\\
	&=&-\frac{\eps}{2}u(u_x+\eps\frac{3}{2}uu_x+\mu a u_{xxt})
	-\mu(a+\frac{1}{3})\partial_{xt}^2(u_x+\eps\frac{3}{2}(uu_x))\\
	& &	+\eps\mu\big[(b+\frac{1}{3})uu_{xx}+(c+\frac{1}{2})
	u_x^2\big]_x
	+O(\mu^2)\\
	&=&-\eps\frac{1}{2}uu_x-\eps^2\frac{3}{4}u^2u_x-\mu(a+\frac{1}{3})u_{xxt}\\
	& &+\eps\mu\big[(2a+b+\frac{5}{6})uu_{xx}+
	(\frac{5}{4}a+c+1)u_x^2\big]_x+O(\mu^2);
\end{eqnarray*}
similarly, one gets
$$
	\eps^2[vu]_x=\eps^2\frac{3}{4}u^2u_x-
	\eps\mu(a+\frac{1}{3})\big[uu_{xx}\big]_x+O(\mu^2),
$$
so that (\ref{eq6}) is equivalent to
$$
	u_t+u_x+\eps\frac{3}{2}uu_x-\mu\eps(a+\frac{1}{3})u_{xxt}
	=\eps\mu\big[
	-(a+b+\frac{1}{2})uu_{xx}-(\frac{5}{4}a+c+1)u_x^2
	\big]+O(\mu^2).
$$
Equating the coefficients of this equation with those of (\ref{eq4})
shows that the first equation of (\ref{eq1}) is also satisfied at 
order $O(\mu^2)$ if the following relations hold:
$$
	a=-\frac{1}{6},\qquad
	b=-\frac{1}{6},\quad
	c=-\frac{19}{48},
$$
and the conditions given in the statement of the proposition on 
$\alpha$, $\beta$, $\gamma$ and $\delta$ follows from the expressions
of $a$, $b$ and $c$ given after equation (\ref{eq4}).
\end{proof}

As said in Remark \ref{rm1}, none of the equations of the one parameter
family considered in Proposition \ref{prop1} is completely integrable.
This is the reason why we now want to derive a wider class of equations
of the form (\ref{CHfamily}) -- and whose solution can still be used as the
basis of an approximate solution of the Green-Naghdi equations
(\ref{eq1}) (and thus of the water waves problem). 
We can generalize Proposition \ref{prop1} by replacing the
vertically averaged velocity $u$ given by (\ref{averaged}) by
the horizontal velocity $u^\theta$ ($\theta\in [0,1]$) evaluated at the level 
line $\theta$ of the fluid domain:
$$
	u^\theta(x)=\dx\Phi_{\vert_{z= (1+\eps\zeta)\theta-1}},
$$
so that $\theta=0$ and $\theta=1$
correspond to the bottom and surface respectively. The introduction of
$\theta$ allows us to derive an approximation consistent with (\ref{eq1})
built on a two-parameter family of equations of the form (\ref{CHfamily}).
\begin{prop}\label{prop2}
	Let $p\in\R$, $\theta\in [0,1]$, and write 
	$\lambda=\frac{1}{2}(\theta^2-\frac{1}{3})$. Asssuming that 
	$$
	\alpha=p+\lambda,
	\quad \beta=p-\frac{1}{6}+\lambda,\quad
	\gamma=-\frac{3}{2}p-\frac{1}{6}-\frac{3}{2}\lambda,\quad
	\delta=-\frac{9}{2}p -\frac{23}{24}-\frac{3}{2}\lambda,
	$$
	there exists $D>0$ such that:
	\begin{itemize}
	\item For all $s\geq 0$ and $T>0$,
	\item  For all bounded family 
	$(u^{\eps,\mu,\theta})_{(\eps,\mu)\in{\mathcal P}}\in 
	C([0,\frac{T}{\eps}];H^{s+D}(\R))$ solving (\ref{CHfamily}),
	\end{itemize}
	the family $(u^{\eps,\mu},\zeta^{\eps,\mu})_{(\eps,\mu)
	\in {\mathcal P}}$,
	with (ommiting the indexes $\eps,\mu$),
	\begin{eqnarray}\label{formulezeta}
	u&=&u^\theta+\mu\lambda u^\theta_{xx}
	+2\mu\eps\lambda u^\theta u^\theta_{xx},\\
	\label{formulebis}
	\zeta&:=&u+\frac{\eps}{4}u^2+\mu\frac{1}{6}u_{xt}
	-\eps\mu\big[\frac{1}{6}uu_{xx}+\frac{5}{48}
	u_x^2\big],
	\end{eqnarray}
	is consistent (of order $s$ and on $[0,\frac{T}{\eps}]$) with 
	the Green-Naghdi equations (\ref{eq1}).
\end{prop}
\begin{rema}\label{rm2}
{\bf i.} The one parameter family of equations (\ref{CHfamily}) of Proposition 
\ref{prop1} corresponds
to the particular case $\theta^2=1/3$ (or $\lambda=0$).\\
{\bf ii.} There exists only one set of coefficients such that
$\delta=2\gamma$ and $\beta=-2\gamma$ (corresponding to $p=-\frac{1}{3}$
and $\theta^2=\frac{1}{2}$, and thus $\alpha=-\frac{1}{4}$,
$\beta=-\frac{5}{12}<0$, $\gamma=\frac{5}{24}$, $\delta=\frac{5}{12}$). 
The corresponding equation is therefore a Camassa-Holm equation
in the sense of Definition \ref{defiCH}:
\begin{equation}\label{equationCH}
	u_t+u_x+\frac{3}{2}\eps u u_x
	-\mu (\frac{1}{4} u_{xxx}+\frac{5}{12} u_{xxt})
	=\frac{5}{24}\eps\mu(u u_{xxx}+2 u_x u_{xx}).
\end{equation}
{\bf iii.} There exists only one set of coefficients such that 
$\delta=3\gamma$ and $\beta=-\frac{8}{3}\gamma$ (obtained with $\theta^2=\frac{23}{36}$ and $p=-\frac{77}{216}$, and thus $\alpha=-\frac{11}{54}$,
$\beta=-\frac{10}{27}$, $\gamma=\frac{5}{36}$, $\delta=\frac{5}{12}$).
The corresponding equation is therefore a Degasperis-Procesi equation
in the sense of Definition \ref{defiDP}:
$$
	u_t+u_x+\frac{3}{2}\eps u u_x
	-\mu (\frac{11}{54} u_{xxx}+\frac{10}{27} u_{xxt})
	=\frac{5}{36}\eps\mu(u u_{xxx}+ 3 u_x u_{xx}).
$$
\end{rema}
\begin{proof}
From the proof of Prop. 3.8 of \cite{AL}, one has the following expression
(neglecting terms of order $O(\mu^2)$):
\begin{eqnarray*}
	u&=&\dx \psi+\mu \big((1+\eps\zeta)\dx\zeta\dx\psi+\frac{(1+\eps\zeta)^2}{3}\dx^3\psi\big),\\
	u^\theta&=&\dx\psi+ \mu \big((1+\eps\zeta)\dx\zeta\dx\psi+\frac{(1+\eps\zeta)^2}{2}(1-\theta^2)\dx^3\psi\big),
\end{eqnarray*}
where $\psi$ denotes the trace of the velocity potential at the surface.
It follows from these formulas that
\begin{eqnarray*}
	u&=&u^\theta+\mu\frac{(1+\eps\zeta)^2}{2}(\theta^2-\frac{1}{3})u^\theta_{xx}
	+O(\mu^2)\\
	&=&u^\theta+\mu\frac{1}{2}(\theta^2-\frac{1}{3})u^\theta_{xx}
	+\mu\eps(\theta^2-\frac{1}{3})u^\theta u^\theta_{xx}
	+O(\mu^2),
\end{eqnarray*}
where we used $\zeta=u^\theta+O(\eps)$ for the last equality.\\
This formula, together with Proposition \ref{prop1}, easily yields the
result.
\end{proof}

\subsection{Equations on the surface elevation}\label{sectsurf}

Proceeding exactly as in the proof of Proposition \ref{prop1},
one can prove that the family of equations 
\begin{equation}\label{CHfamilyzeta}
	\zeta_t+\zeta_x+\frac{3}{2}\eps \zeta \zeta_x
	-\frac{3}{8}\eps^2\zeta^2\zeta_x
	+\frac{3}{16}\eps^3\zeta^3\zeta_x
	+\mu (\alpha \zeta_{xxx}+\beta \zeta_{xxt})
	=\eps\mu(\gamma \zeta \zeta_{xxx}+\delta \zeta_x\zeta_{xx})
\end{equation}
for the evolution of the surface elevation can be used to construct an approximate
solution consistent with the Green-Naghdi equations:
\begin{prop} \label{prop1zeta}
	Let $q\in\R$ and asssume that 
	$$
	\alpha=q,\quad \beta=q-\frac{1}{6},\quad
	\gamma=-\frac{3}{2}q-\frac{1}{6},\quad
	\delta=-\frac{9}{2}q-\frac{5}{24}.
	$$
	Then there exists $D>0$ such that:
	\begin{itemize}
	\item For all $s\geq 0$ and $T>0$,
	\item  For all bounded family 
	$(\zeta^{\eps,\mu})_{(\eps,\mu)\in{\mathcal P}}\in 
	C([0,\frac{T}{\eps}];H^{s+D}(\R))$ solving (\ref{CHfamilyzeta}),
	\end{itemize}
	the family $(\zeta^{\eps,\mu},u^{\eps,\mu})_{(\eps,\mu)
	\in {\mathcal P}}$,
	with (omitting the indexes $\eps,\mu$)
	$$
	u:=\zeta+\frac{1}{h}\Big(-\frac{\eps}{4}\zeta^2
	-\frac{\eps^2}{8}\zeta^3
	+\frac{\eps^3}{64}\zeta^4
	-\mu\frac{1}{6}\zeta_{xt}
	+\eps\mu\big[\frac{1}{6}\zeta\zeta_{xx}+\frac{1}{48}
	\zeta_x^2\big]\Big),
	$$
	is consistent (of order $s$ and on $[0,\frac{T}{\eps}]$) with 
	the Green-Naghdi equations (\ref{eq1}).
\end{prop}
\begin{rema}\label{remzeta}
	Choosing $q=1/12$, the equation (\ref{CHfamilyzeta}) reads
	\begin{eqnarray}\label{eqs3}
	&&\zeta_t+\zeta_x+\frac{3}{2}\eps \zeta \zeta_x
	-\frac{3}{8}\eps^2\zeta^2\zeta_x
	+\frac{3}{16}\eps^3\zeta^3\zeta_x
	+\frac{\mu}{12} (\zeta_{xxx}- \zeta_{xxt})\\
	&&\hskip 5cm =-\frac{7}{24}\eps\mu(\zeta \zeta_{xxx}+ 2\zeta_x\zeta_{xx}).\nonumber
	\end{eqnarray}
While for any $q \in \R$,  (\ref{CHfamilyzeta}) is an equation for the evolution of the free surface $\zeta$, 
and all these equations have  
the same order of accuracy $O(\eps^4,\mu^2)$, it is more advantageous to use (\ref{eqs3}) since it presents 
better structural properties that we take advantage of in \S \ref{wb}. The ratio $2:1$ between 
the coefficients of $\zeta_x\zeta_{xx}$ and $\zeta\zeta_{xxx}$ is crucial in our considerations.	
\end{rema}

\section{Mathematical analysis of the models and rigorous justification}

\subsection{Large time well-posedness of the unidirectional equations
(\ref{CHfamily}) and (\ref{CHfamilyzeta})}

We prove here the well posedness of the general class of equations 
\begin{eqnarray}\label{CHgeneral}
&u_t+u_x+\frac{3}{2}\eps u u_x+\eps^2\iota u^2u_x+\eps^3\kappa u^3u_x +
	\mu (\alpha u_{xxx}+\beta u_{xxt})\\
&\qquad =\eps\mu(\gamma u u_{xxx}+\delta u_x u_{xx}),\nonumber
\end{eqnarray}
with $\iota,\,\kappa \in \R$; in particular, (\ref{CHgeneral}) coincides
with (\ref{CHfamily}) and (\ref{CHfamilyzeta}) if one takes
$\iota=\kappa=0$ and $\iota=-\frac{3}{8}$,
$\kappa=\frac{3}{16}$ respectively. That is, we solve the initial
value problem
\begin{equation}\label{ivpCH}
	\left\vert
	\begin{array}{l}
	u_t+u_x+\frac{3}{2}\eps u u_x+\eps^2\iota u^2u_x+\eps^3\kappa u^3u_x+
	\mu (\alpha u_{xxx}+\beta u_{xxt})\\
	\hskip 3cm =\eps\mu(\gamma u u_{xxx}+\delta u_x u_{xx}),\\
	u_{\vert_{t=0}}=u^0
	\end{array}\right.
\end{equation}
on a time scale
$O(1/\eps)$, and under the condition $\beta<0$. In order to state the
result, we need to define the spaces $X^s$ as
$$
	\forall s\geq 0, \qquad
	X^{s+1}=H^{s+1}(\R)\mbox{ endowed with the norm }\vert f\vert_{X^{s+1}}^2
	=\vert f\vert_{H^s}^2+\mu\vert \dx f\vert_{H^s}^2,
$$
and we also recall that the set ${\mathcal P}$ is defined in (\ref{cond1}).
\begin{prop}\label{prop3}
	Assume that $\beta<0$ and let $\mu_0>0$, $M>0$, 
	$s>\frac{3}{2}$ and $u^0\in H^{s+1}(\R)$.
	Then, there exists $T>0$ and a unique family of solutions 
	$(u_{\eps,\mu})_{(\eps,\mu)\in{\mathcal P}}$ to (\ref{ivpCH})
	bounded in $C([0,\frac{T}{\eps}];X^{s+1}(\R))\cap
	C^1([0,\frac{T}{\eps}];X^{s}(\R))$.
\end{prop}
\begin{proof}
For all $v$ smooth enough, let us define the ``linearized'' operator
${\mathcal L}(v,\partial)$ as
\begin{eqnarray*}
	{\mathcal L}(v,\partial)&=&(1+\mu \beta\dx^2)\dt
	+\dx+\mu\alpha\dx^3+\frac{3}{2}\eps v\dx+\eps^2\iota v^2\dx
	+\eps^3\kappa v^3\dx\\
	& &
	-\eps\mu\gamma v\dx^3
	-\eps\mu\delta\big(\frac{1}{2}v_x\dx^2+\frac{1}{2}v_{xx}\dx\big),
\end{eqnarray*}
In order to construct a solution to (\ref{CHgeneral}) by an iterative scheme,
we are led to study the initial value problem
\begin{equation}\label{ivp}
	\left\lbrace
	\begin{array}{l}
	{\mathcal L}(v,\partial)u=\eps f,\\
	u_{\vert_{t=0}}=u^0.
	\end{array}\right.
\end{equation}
If $v$ is smooth enough, 
it is completely standard to check that for all $s\geq 0$,
$f \in L^1_{loc}(\R^+_t;H^s(\R_x))$ and $u^0\in H^s(\R)$, there
exists a unique solution $u\in C(\R^+;H^{s+1}(\R))$ to (\ref{ivp})
(recall that $\beta<0$). We take for granted the existence of a
solution to (\ref{ivp}) and establish some precise energy estimates
on the solution. In order to do so, let us define the ``energy'' norm
$$
	\forall s\geq 0,\qquad
	E^s(u)^2=\vert u\vert_{H^s}^2-\mu\beta\vert \dx u\vert_{H^s}^2.
$$
Differentiating $\frac{1}{2}e^{-\eps\lambda t}E^s(u)$
with respect to time, one gets, using the equation (\ref{ivp})
and integrating by parts,
\begin{eqnarray*}
	\frac{1}{2}e^{\eps\lambda t}\dt (e^{-\eps\lambda t}E^s(u)^2)
	&=&-\frac{\eps\lambda}{2} E^s(u)^2
	+ \eps (\Lambda^s f,\Lambda^su)
	-\eps (\Lambda^s(V\dx u),\Lambda^su)\\
	& &+\eps\mu\gamma(\Lambda^s(v\dx^3 u),\Lambda^s u)
	-\eps\mu\frac{\delta}{2}
	(\Lambda^s(v_x\dx u),\Lambda^s \dx u),
\end{eqnarray*}
with $V=\frac{3}{2}v+\eps\iota v^2+\eps^2\kappa v^3$.\\
Since for all constant coefficient skewsymmetric differential 
polynomial $P$ (that is, $P^*=-P$), and all $h$ smooth enough, one has
$$
	(\Lambda^s (h P u),\Lambda^s u)=
	([\Lambda^s,h]Pu,\Lambda^s u)
	-\frac{1}{2}([P,h]\Lambda ^su,\Lambda^s u),
$$
we deduce (applying this identity with $P=\dx$ and $P=\dx^3$),
\begin{eqnarray*}
	\lefteqn{\frac{1}{2}e^{\eps\lambda t}\dt(e^{-\eps\lambda t} E^s(u)^2)
	=-\frac{\eps\lambda}{2} E^s(u)^2
	-\eps \big([\Lambda^s,V]\dx u,\Lambda^su)
	+\frac{\eps}{2}((\dx V)\Lambda^s u,\Lambda^s u\big)}\\
	& &-\eps\mu\gamma\big([\Lambda^s,v]\dx^2 u-\frac{3}{2}v_x\Lambda^s\dx u-v_{xx}\Lambda^s u,\Lambda^s \dx u\big)
	-\eps\mu\gamma\big([\Lambda^s,v_x]\dx^2 u,\Lambda^su\big)\\
	& &-\eps\mu\frac{\delta}{2}
	(\Lambda^s(v_x\dx u),\Lambda^s \dx u)
	+ \eps (\Lambda^s f,\Lambda^su),
\end{eqnarray*}
Here we also used the identities
$$[\Lambda^s, v]\dx^3 u=\dx\big([\Lambda^s,v]\dx^2 u\big)-[\Lambda^s,v_x]\dx^2 u
$$
and 
$$
	\frac{1}{2}(v_{xxx}\Lambda^s u,\Lambda^s u)=-(v_{xx}\Lambda^s u,\Lambda^s u_x).
$$
Since $\vert u\vert_{H^s}\leq E^s(u)$ and $\sqrt{\mu}\vert\dx u\vert_{H^s}\leq 
\frac{1}{\sqrt{-\beta}}E^s(u)$,
one gets directly by the Cauchy-Schwartz inequality,
\begin{eqnarray*}
	e^{\eps\lambda t}\dt(e^{-\eps\lambda t} E^s(u)^2)
	&\leq&
	\eps C(\mu_0,\frac{1}{\beta},\gamma,\delta)
	(A(u,v)E^s(u)+B(v)E^s(u)^2)\\
	& & -\eps\lambda E^s(u)^2+2\vert f\vert_{H^s}E^s(u),
\end{eqnarray*}
with
\begin{eqnarray*}
	A(u,v)&=&\vert [\Lambda^s,V]\dx u\vert_2
	+\vert [\Lambda^s,v]\dx (\sqrt{\mu}\dx u)\vert_2
	+\vert [\Lambda^s,\sqrt{\mu}v_x]\dx (\sqrt{\mu}\dx u)\vert_2\\
	& &+\sqrt{\mu}\vert v_x\dx u\vert_{H^s},\\
	B(v)&=&\vert \dx V\vert_\infty+\vert v_x\vert_\infty+\vert \dx(\sqrt{\mu}\dx v)\vert_\infty.
\end{eqnarray*}
Recalling that for all $s> 3/2$, and all $F,U$ smooth enough, one has
$$
	\vert [\Lambda^s,F]U\vert_2\leq \cst \vert F\vert_{H^s}\vert U\vert_{H^{s-1}},
$$
it is easy to check that one gets
 $A(u,v)\leq C(\mu_0,M,\frac{1}{\beta},\iota,\kappa, E^s(v))E^s(u)$ 
and 
$B(v)\leq  C(\mu_0,M,\frac{1}{\beta},\iota, \kappa,E^s(v))$. Therefore, 
we obtain
$$
	e^{\eps\lambda t}\dt (e^{-\eps\lambda t}E^s(u)^2) 
	\leq  \big(C(\mu_0,M,\frac{1}{\beta},\gamma,\delta,\iota,\kappa,
	E^s(v))-\lambda\big)E^s(u)^2+2\eps E^s(f)E^s(u).
$$
Taking $\lambda=\lambda_T$ large enough (how large depending on $C(\mu_0,M,\frac{1}{\beta},\gamma,\delta,\iota,\kappa,E^s(v))$) to have the first term of the right hand side
negative for all $t\in [0,\frac{T}{\eps}]$, one deduces
$$
	\forall t\in [0,\frac{T}{\eps}],\qquad
	\dt (e^{-\eps\lambda_T t}E^s(u)^2)\leq 
	2\eps e^{-\eps\lambda_T t}E^s(f)E^s(u).
$$
Integrating this differential inequality yields therefore
$$
	\forall t\in [0,\frac{T}{\eps}],\qquad
	E^s(u)(t)\leq e^{\eps\lambda_T t}E^0(u^0)
	+2\eps \int_0^t e^{\lambda_T (t-t')}
	E^s(f(t'))dt'.
$$
Thanks to this energy estimate, one can conclude classically
(see e.g. \cite{AG})
to the existence of 
$$
	T=T(\mu_0,M,\vert u^0\vert_{X^{s+1}_{\mu_0}},\frac{1}{\beta},\gamma,\delta,\iota,\kappa)>0,
$$ 
and of
a unique solution $u\in C([0,\frac{T}{\eps}];X^{s+1}(\R^d))$ to
(\ref{ivpCH}) as a limit of the iterative scheme
$$
	u_0=u^0,\quad\mbox{ and }\quad
	\forall n\in\N, \quad
	\left\lbrace\begin{array}{l}
	{\mathcal L}(u^n,\partial)u^{n+1}=0,\\
	u^{n+1}_{\vert_{t=0}}=u^0.
		    \end{array}\right.
$$
Since $u$ solves (\ref{CHgeneral}), we have ${\mathcal L}(u,\partial u)u=0$
and therefore 
$$
	(\Lambda^{s-1}(1+\mu\beta\dx^2)\dt u,\Lambda^{s-1}\dt u)
	=-\eps(\Lambda^{s-1}{\mathcal M}(u,\partial)u,\Lambda^{s-1}\dt u),
$$
with ${\mathcal M}(u,\partial)={\mathcal L}(u,\partial)-(1+\mu\beta\dx^2)\dt$.
Proceeding as above, one gets
$$
	E^{s-1}(\dt u)\leq
	C(\mu_0,M,\vert u^0\vert_{X^{s+1}_{\mu_0}},\frac{1}{\beta},\gamma,\delta,\iota,\kappa, E^s(u)),
$$
and it follows that the family of solution is also bounded in 
$C^1([0,\frac{T}{\eps}];X^s)$.
\end{proof}

\subsection{Rigorous justification of the unidirectional approximations (\ref{CHfamily})}

In Proposition \ref{prop2}, we constructed a family 
$(u^{\eps,\mu},\zeta^{\eps,\mu})$ \emph{consistent} with the
Green-Naghdi equations in the sense of Definition \ref{defi1}.
A consequence of the following theorem is a stronger result: this family
provides a good approximation of the exact solutions 
$(\underline{u}^{\eps,\mu},\underline{\zeta}^{\eps,\mu})$ of the
Green-Naghdi equations with same initial data in the sense that
$(\underline{u}^{\eps,\mu},\underline{\zeta}^{\eps,\mu})=(u^{\eps,\mu},\zeta^{\eps,\mu})+O(\mu^2 t)$ for times $O(1/\eps)$.
\begin{theo}
	Let $\mu_0>0$, $M>0$, $T>0$ and ${\mathcal P}$ 
	be as defined in (\ref{cond1}).
	Let also $p\in\R$, $\theta\in [0,1]$, and
	$\alpha$, $\beta$, $\gamma$ and $\delta$ be as in Proposition 
	\ref{prop2}.\\
	If $\beta<0$ then
	there exists $D>0$ and $T>0$ such that for all 
	$u^{0}\in H^{s+D+1}(\R)$:
	\begin{itemize}
	\item There is a unique family 
	$(u^{\eps,\mu},\zeta^{\eps,\mu})_{(\eps,\mu)\in{\mathcal P}}\in 
	C([0,\frac{T}{\eps}];H^{s+D}(\R)^2)$ given by 
	the resolution of (\ref{CHfamily}) 
	with initial condition $u^{0}$ and formulas (\ref{formulezeta})-(\ref{formulebis});
	\item  There is a unique family
	$(
	\underline{u}^{\eps,\mu},\underline{\zeta}^{\eps,\mu})_{(\eps,\mu)\in{\mathcal P}}\in 
	C([0,\frac{T}{\eps}];H^{s+D}(\R)^2)$ solving the Green-Naghdi equations
	(\ref{eq1}) with initial condition 
	$(u^{\eps,\mu},\zeta^{\eps,\mu})_{\vert_{t=0}}$.
	\end{itemize}
	Moreover, one has for all $(\eps,\mu)\in {\mathcal P}$,
	$$
	\forall t\in [0,\frac{T}{\eps}],\qquad
	\vert \underline{u}^{\eps,\mu}-u^{\eps,\mu}\vert_{L^\infty([0,t]\times \R)}+\vert \underline{\zeta}^{\eps,\mu}-\zeta^{\eps,\mu}\vert_{L^\infty([0,t]\times \R)}\leq \cst \mu^2 t.
	$$
\end{theo}
\begin{rema}
	It is known (see \cite{AL}) that the Green-Naghdi equations give,
	under the scaling (\ref{scalingCH}), a correct approximation of 
	the exact solutions of the full water waves equations (with a precision
	$O(\mu^2 t)$ and over a time scale $O(1/\eps)$). It follows that
	that the unidirectional approximation discussed above
	approximates the solution of the water waves equations
	with the same accuracy.
\end{rema}
\begin{rema}
	We used the unidirectional equations derived on the velocity
	as the basis for the approximation justified in the theorem.
	One could of course use instead the unidirectional approximation
	(\ref{CHfamilyzeta}) derived on the surface elevation.
\end{rema}
\begin{proof}
The first point of the theorem is a direct consequence of Proposition \ref{prop3}. Thanks to Proposition \ref{prop2}, we now that $(u^{\eps,\mu},\zeta^{\eps,\mu})_{\eps,\mu}$ is consistent with the Green-Naghdi equations (\ref{eq1}),
so that the second point of the theorem and the error estimate follow
at once from the well-posedness and stability of the Green-Naghdi
equations (see Th. 3 of \cite{AL2}-- note that instead of using this
general result which holds for two dimensional surfaces and nonflat bottoms,
one could easily adapt the simpler and more precise results of \cite{Li}
to the present scaling).
\end{proof}

\subsection{Wave breaking}\label{wb}

For the Camassa-Holm family of equations (\ref{CHfamily}) for the velocity it is known 
(see \cite{Co}) that singularities can develop in finite time for a smooth initial data only 
in the form of wave breaking. We will show now that this form of blow-up is also a 
feature of the equation (\ref{eqs3}) for the free surface. More precisely, if a smooth 
initial profile fails to produce a wave that exists for all subsequent times, then we 
encounter wave breaking in the form of surging (and not plunging, as would be the 
case if (\ref{CHfamily}) were the equation for the evolution of the free surface). 

Our first result describes the precise blow-up pattern for the equation 
(\ref{eqs3}) for the free surface.

\begin{prop}\label{blowCH1}
Let $\zeta_0 \in H^3(\R)$. If the maximal existence time $T>0$ of the solution 
of (\ref{eqs3}) with initial profile $\zeta(0,\cdot)=\zeta_0$ is finite, $T< \infty$, then 
the solution $\zeta \in C^1([0,T);H^2(\R)) \cap C([0,T);H^(\R))$ is such that
\begin{equation}\label{CHb1}
\sup_{t \in [0,T),\,x \in \R}\{ |\zeta(t,x)|\} < \infty
\end{equation}
and 
\begin{equation}\label{CHb2}
\sup_{x \in \R}\,\{ \zeta_x(t,x)\} \uparrow \infty\quad as \quad t \uparrow T.
\end{equation}
\end{prop}

\begin{proof}
In view of Proposition \ref{prop3}, given $\zeta_0 \in H^3(\R)$, the maximal 
existence time of the solution $\zeta(t)$ to (\ref{eqs3}) with initial data $\zeta(0)=\zeta_0$ 
is finite if and only if $\vert \zeta(t)\vert_{H^3(\R)}$ blows-up in finite time. Thus if (\ref{CHb2}) holds for some 
finite $T>0$, then the maximal existence time is finite. To complete the proof it suffices to show 
that 

(i) the solution $\zeta(t)$ given by Proposition \ref{prop3} remains uniformly bounded as long as it is defined;

\noindent
and 

(ii) if we can find some $M=M(\zeta_0)>0$ such that 
\begin{equation}\label{CHb3}
\zeta_x(t,x) \le M,\qquad x \in \R,
\end{equation}
as long as the solution is defined, then $\vert \zeta(t)\vert_{H^3(\R)}$ stays bounded on bounded time-intervals. 

Item (i) follows at once from the imbedding $L^\infty(\R) \subset H^1(\R)$ since multiplying (\ref{eqs3}) 
by $\zeta$ and integrating on $\R$ yields
\begin{equation}\label{CHb4}
\partial_t\Big(\int_\R [\zeta^2+\frac{1}{12}\,\mu\int_\R \zeta_x^2]\,dx\Big)=0.
\end{equation}

To prove item (ii), notice that multiplication of (\ref{eqs3}) by $\zeta_{xxxx}$ 
yields
\begin{eqnarray}\label{CHb5}
&&\partial_t\Big( \int_\R [\zeta_{xx}^2+\frac{1}{12}\,\mu\int_\R \zeta_{xxx}^2]\,dx\Big) =15\,\eps\int_{\R}\zeta\zeta_{xx}\zeta_{xxx}\,d-\frac{15}{4}
\,\eps^2 \int_{\R} \zeta^2\zeta_{xx}\zeta_{xxx}\,dx \\
&&\qquad
 +\frac{9}{16}\eps^3 \int_\R \zeta_x^5\,dx
+\frac{15}{8}\,\eps^3\int_\R\zeta^3\zeta_{xx}\zeta_{xxx}\,dx+\frac{7}{4}\,\mu\eps\int_\R \zeta_x\zeta_{xxx}^2\,dx.\nonumber
\end{eqnarray}
after performing several integrations by parts. 
If (\ref{CHb3}) holds, let in accordance with (\ref{CHb4}) the constant $M_0>0$ be such that 
$$|\zeta(t,x)| \le M_0,\qquad x \in \R,$$
for as long as the solution exists. Using the Cauchy-Schwartz inequality 
as well as the fact that $\mu \le 1$, we infer from (\ref{CHb4}) and (\ref{CHb5}) that
$$\partial_t E(t)\le \Big(\frac{90\eps}{\mu}\,M_0+\frac{45\eps^2}{2\mu}\,M_0^2+\frac{27\eps^3}{4\mu}\,M^3 
+\frac{45\eps^3}{4\mu}\,M_0^3+\,21\,\eps M  \Big)\, E(t),$$
where
$$E(t)=\int_\R [\zeta^2+ \frac{1}{12}\,\mu\zeta_x^2+ \zeta_{xx}^2 +\frac{1}{12}\,\mu
 \zeta_{xxx}^2]\,dx.$$
An application of Gronwall's inequality enables us to conclude.
\end{proof}

Our next aim is to show that there are solutions to (\ref{eqs3}) that blow-up in finite time as 
surging breakers, that is, following the pattern given in Proposition \ref{blowCH1}. 
We will prove this by analyzing the equation that describes the evolution of 
\begin{equation}\label{CHb6}
M(t)=\sup_{x \in \R}\,\{\zeta_x(t,x)\}.
\end{equation}
For the degree of smoothness of the solution $\zeta(t)$ given by Proposition \ref{prop3}, we know 
that $M(t)$ is locally Lipschitz and
\begin{equation}\label{CHb7}
\frac{d}{dt}\,M(t)=\zeta_{tx}(t,\xi(t))\quad\hbox{for a.e.}\quad t,
\end{equation}
where $\xi(t)$ is any point where $M(t)=\zeta_x(t,\xi(t))$ cf. \cite{CE}. For further use, let us also note that 
\begin{equation}\label{CHb8}
(1-\frac{1}{12}\,\mu\,\partial_x^2)^{-1}f=P \ast f,\qquad f \in L^2(\R),
\end{equation}
where
$$P(x)=\sqrt{\frac{3}{\mu}}\ e^{-\displaystyle\,2\,\sqrt{\frac{3}{\mu}}\,|x|},\qquad x \in \R,$$
with
\begin{equation}\label{CHb9}
\| P\|_{L^\infty}=\sqrt{\frac{3}{\mu}},\quad \|P\|_{L^1}=1,\quad \|P\|_{L^2}=\Big(\frac{3}{4\mu}\Big)^{\frac{1}{4}},
\end{equation}
and
\begin{equation}\label{CHb10}
\quad \|P_x\|_{L^\infty}=\frac{6}{\mu},
\quad \|P_x\|_{L^1}=2\,\sqrt{\frac{3}{\mu}},\quad \|P_x\|_{L^2}=
\sqrt{2}\,\Big(\frac{3}{\mu}\Big)^{\frac{3}{4}} \le 4\,\mu^{-\frac{3}{4}}.
\end{equation}
Applying $(1-\frac{1}{12}\,\mu\,\partial_x^2)^{-1}$ to (\ref{eqs3}), we obtain the equation 
\begin{eqnarray*}
&&\zeta_t+P_x \ast \zeta +\frac{3}{4}\,\eps\,P_x \ast \zeta^2 - \frac{1}{8}\,\eps^2\,P_x \ast \zeta^3 
+\frac{3}{64}\,\eps^3\,P_x \ast \zeta^4\\
&&\qquad + \frac{1}{12}\,\mu\,\partial_x^3\,P \ast \zeta 
 =-\frac{7}{24}\,\mu\eps\,P_x \ast \zeta_x^2 - \frac{7}{24}\,\mu\eps \, P \ast (\zeta\zeta_{xxx}).
\end{eqnarray*}
Differentiating this equation with respect to the spatial variable, we obtain
\begin{eqnarray*}
&&\zeta_{tx}+\partial_x^2 P \ast \zeta +\frac{3}{4}\,\eps\,\partial_x^2 P\ast \zeta^2 - \frac{1}{8}\,\eps^2\,
\partial_x^2 P \ast \zeta^3 
+\frac{3}{64}\,\eps^3\,\partial_x^2 P \ast \zeta^4 \\
&&\qquad + \frac{1}{12}\,\mu\,\partial_x^4\,P \ast \zeta 
 =-\frac{7}{24}\,\mu\eps\,\partial_x^2 P \ast \zeta_x^2 - \frac{7}{24}\,\mu\eps \, P_x \ast (\zeta\zeta_{xxx}).
\end{eqnarray*}
Since $\zeta\zeta_{xxx}=\partial_x^2(\zeta\zeta_x)-3\zeta_x\zeta_{xx}$ and
\begin{equation}\label{Chb11}
\partial_x^2 \, P\ast f=P_x \ast \zeta_x=\frac{12}{\mu}\, P\ast f - \frac{12}{\mu}\,f,\qquad f \in L^2(\R),
\end{equation}
we deduce that
\begin{eqnarray}\label{CHb12}
&&\zeta_{tx}+2\,P_x \ast \zeta_x-\frac{3}{8}\, 
\eps^2\,P_x\ast (\zeta^2\zeta_x) +\frac{3\eps^3}{16} \,P_x\ast (\zeta^3\zeta_x)\\
&&\qquad =\zeta_{xx} +\frac{7\eps}{4}\,P \ast \zeta_x^2 +\frac{7\eps}{4}\,\zeta_x^2+\frac{7\eps}{2}\,\zeta\zeta_{xx}-\frac{7}{2}\,\eps\,P_x \ast (\zeta\zeta_x).
\nonumber
\end{eqnarray}
We can now prove the following blow-up result.

\begin{prop}\label{blowCH2}
If the initial wave profile $\zeta_0 \in H^3(\R)$ satisfies
\begin{eqnarray*}
\Big|\sup_{x \in \R} \,\{\zeta_0(x)\}\Big|^2 &\ge& \frac{28}{3}\,C_0\,\mu^{-3/4}
+\frac{1}{2}\,\eps\,C_0^{3/2}\,\mu^{-3/4}+\frac{1}{4}\,\eps^2\,C_0^2\,\mu^{-3/4}\\
&&\qquad +\frac{7}{3}\,C_0\,\mu^{-1/2}+\frac{16}{3}\,C_0^{1/2}\,\mu^{-3/4}\,\eps^{-1},
\end{eqnarray*}
where
$$C_0=\int_\R [\zeta_0^2+(\zeta_0')^2]\,dx>0,$$
then wave breaking occurs for the solution of (\ref{eqs3}) in finite time $T=O(\frac{1}{\eps})$.
\end{prop}

\begin{proof}
Notice that
$$\sup_{x \in \R}\,\{|\zeta^2(x)|\} \le \frac{1}{2}\,\int_\R (\zeta^2+\zeta_x^2)\,dx=\frac{C_0}{2}.$$
Therefore, using Young's inequality and the estimates (\ref{CHb9})-(\ref{CHb10}), we obtain that 
\begin{eqnarray*}
&&\| P_x \ast \zeta_x \|_{L^\infty} \le \| P_x \|_{L^2}\| \zeta_x \|_{L^2} \le 4\,\mu^{-3/4}\,C_0^{1/2},\\
&&\| P \ast \zeta_x^2 \|_{L^\infty} \le \| P \|_{L^\infty}\| \zeta_x^2 \|_{L^1} 
\le \| P \|_{L^\infty}\| \zeta_x \|_{L^2}^2 \le 2\,\mu^{-1/2}\,C_0,\\
&&\| P_x \ast (\zeta\zeta_x) \|_{L^\infty} \le \| P_x \|_{L^2}\| \zeta\zeta_x \|_{L^2} \le 
\| P_x \|_{L^2}\| \zeta \|_{L^\infty}\|\zeta_x \|_{L^2}\le 4\,\mu^{-3/4}\,C_0,\\
&&\| P_x \ast (\zeta^2\zeta_x) \|_{L^\infty} \le \| P_x \|_{L^2}\| \zeta^2\zeta_x \|_{L^2} \le 
\| P_x \|_{L^2}\|\zeta\|_{L^\infty}^2\| \|\zeta_x \|_{L^2}\le 2\,\mu^{-3/4}\,C_0^{3/2},\\
&&\| P_x \ast (\zeta^3\zeta_x) \|_{L^\infty} \le \| P_x \|_{L^2}\| \zeta^3\zeta_x \|_{L^2} \le 
\| P_x \|_{L^2}\|\zeta\|_{L^\infty}^3\| \|\zeta_x \|_{L^2}\le 2\, \mu^{-3/4}\,C_0^2.
\end{eqnarray*}
Since (\ref{CHb12}) is at any fixed time an equality in the space of continuous functions, we can 
evaluate both sides at some fixed time $t$ at a point $\xi(t) \in \R$ where $M(t)=\zeta_x(t,\xi(t))$, with 
$M(t)$ defined in (\ref{CHb6}). Since $\zeta_{xx}(t,\xi(t))=0$, from (\ref{CHb7}), (\ref{CHb12}) and 
the previous estimates we derive 
the following differential inequalities for the locally Lipschitz function $M(t)$:
\begin{eqnarray}\label{CHb13}
&&\frac{7}{4}\,\eps\,M^2(t)+\Big(14\,C_0\,\eps+\frac{3}{4}\,\eps^2\,C_0^{3/2}+ 
\frac{3}{8}\,\eps^3\,C_0^2+8C_0^{1/2}\Big)\,\mu^{-3/4} +\frac{7}{2}\,\eps C_0\mu^{-1/2} \\
&&\quad \ge M'(t) \quad\hbox{for a.e.}\quad t,\nonumber
\end{eqnarray}
and
\begin{eqnarray}\label{CHb14}
&&M'(t) \ge \frac{7}{4}\,\eps\,M^2(t)-\Big(14\,C_0\,\eps+\frac{3}{4}\,\eps^2\,C_0^{3/2}+ 
\frac{3}{8}\,\eps^3\,C_0^2+8C_0^{1/2}\Big)\,\mu^{-3/4}\\
&&\hskip 5cm \quad\hbox{for a.e.}\quad t.\nonumber
\end{eqnarray}
Notice that $0 \not \equiv \zeta_0 \in H^3(\R)$ ensures $M(0)>0$. At $t=0$ the right-hand side 
of (\ref{CHb14}) is by our assumption on the initial wave profile larger than $\displaystyle\frac{1}{4}\,\eps\,M^2(0)$. We infer that up to the maximal existence time $T>0$ of 
the solution $\zeta(t)$ of (\ref{eqs3}) the function $M(t)$ must be increasing and, moreover,
$$M'(t) \ge \frac{1}{4}\,\eps\,M^2(t) \quad\hbox{for a.e.}\quad t.$$
Dividing by $M^2(t) \ge M^2(0)>0$ and integrating, we get
$$\frac{1}{M(0)}-\frac{1}{M(t)} \ge \frac{1}{4}\,\eps\,t,\qquad t \in [0,T).$$
Therefore $\lim_{t \uparrow T}\,M(t)=\infty$ and $T \le \displaystyle\frac{4}{\eps\,M(0)}$. 

On the other hand, a similar argumentation applied to (\ref{CHb13}) yields
$$M'(t) \le 4\eps\,M^2(t) \quad\hbox{for a.e.}\quad t,$$
so that
$$\frac{1}{M(t)} \ge \frac{1}{M(0)}-4\eps\, t$$
as long as the solution of (\ref{eqs3}) is defined. Since $\lim_{t \uparrow T}\,M(t)=\infty$ we 
deduce from the previous inequality that $T \ge \displaystyle\frac{1}{4\eps\,M(0)}$. Thus the 
finite maximal existence time $T>0$ is of order $O(\frac{1}{\eps})$.
\end{proof}

\subsection{Numerical computations}\label{numerics}

\begin{figure} [htb]
	\begin{center}
	\includegraphics[width=2.9cm,angle=270]{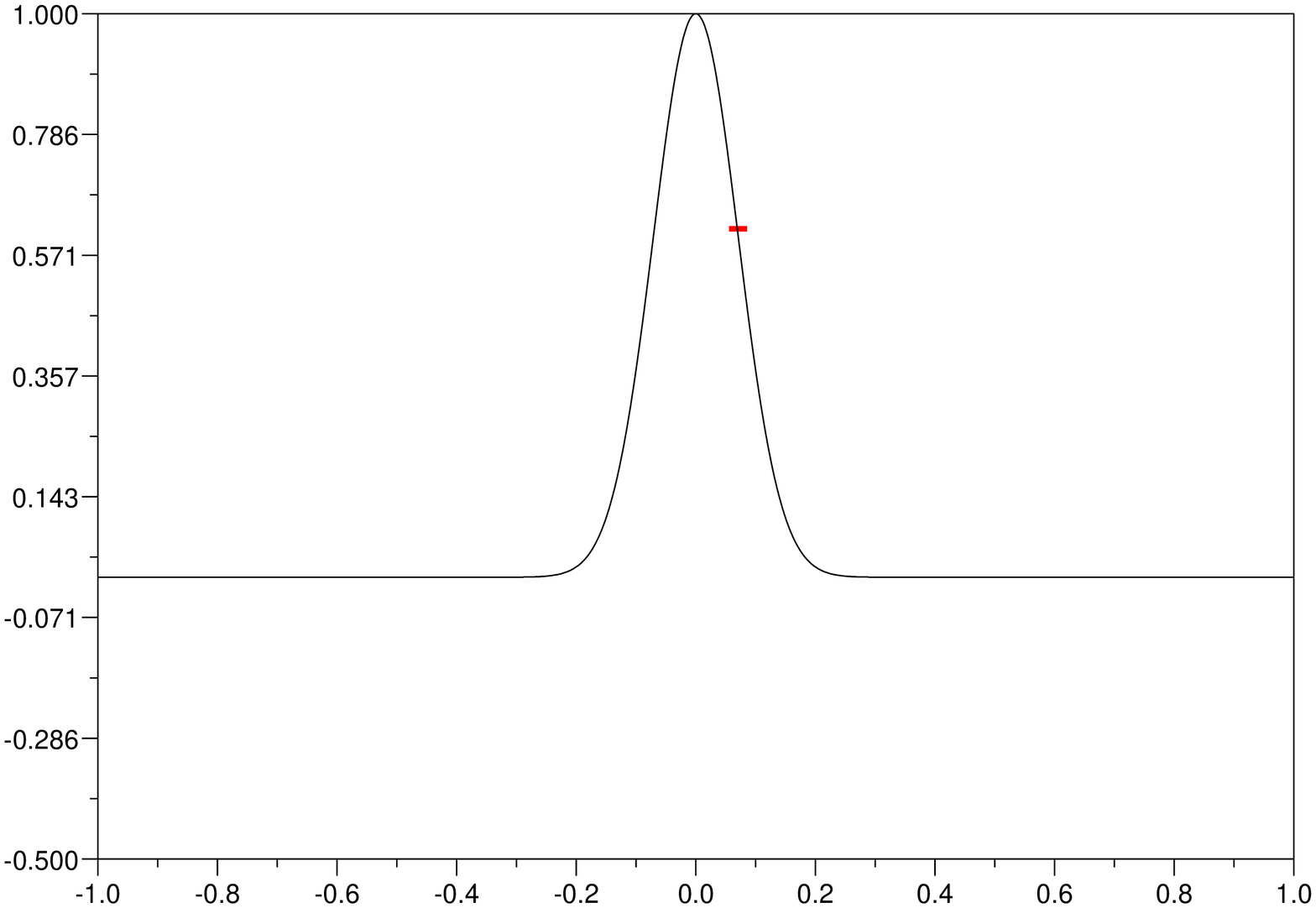}
	\includegraphics[width=2.9cm,angle=270]{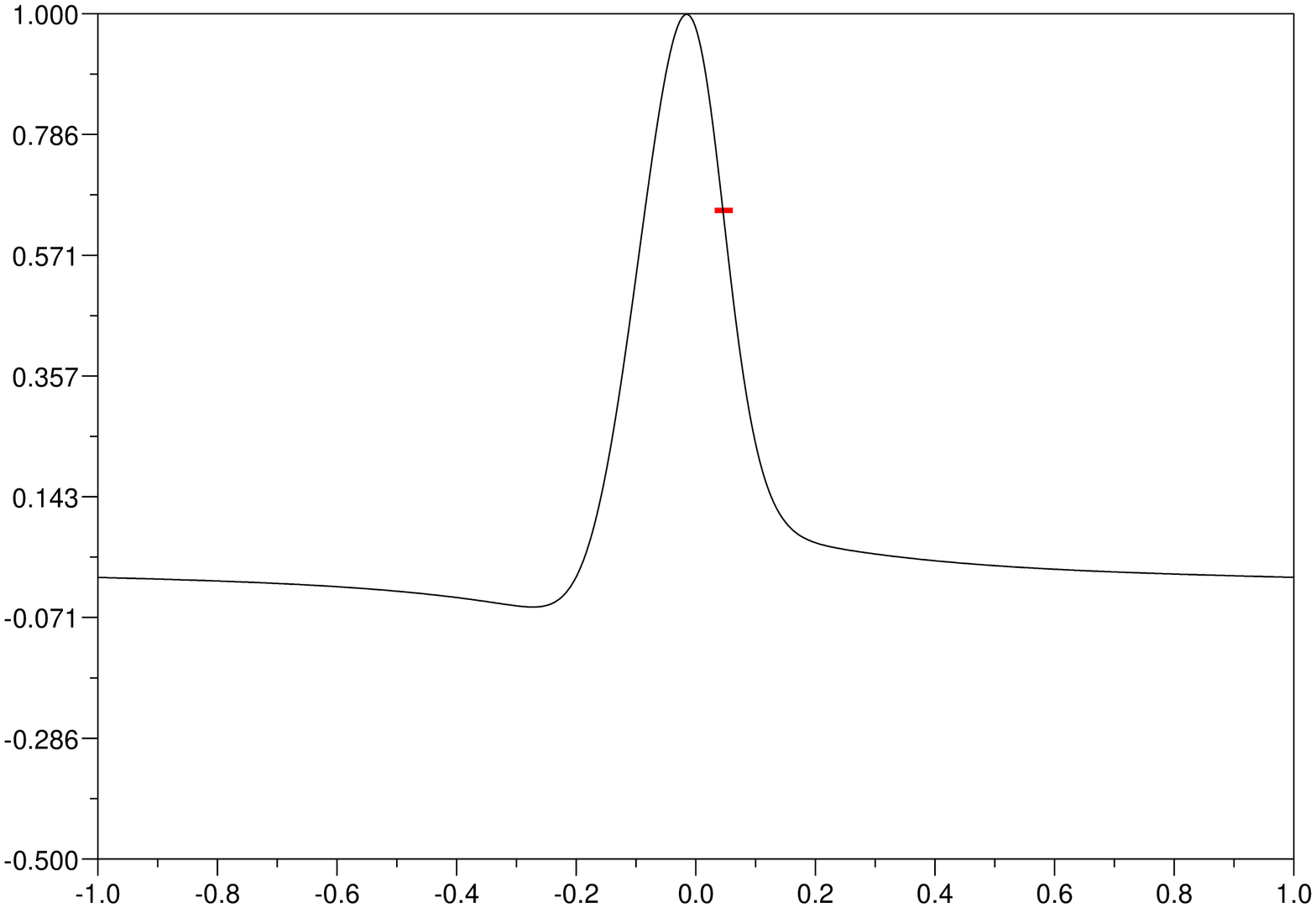}
	\includegraphics[width=2.9cm,angle=270]{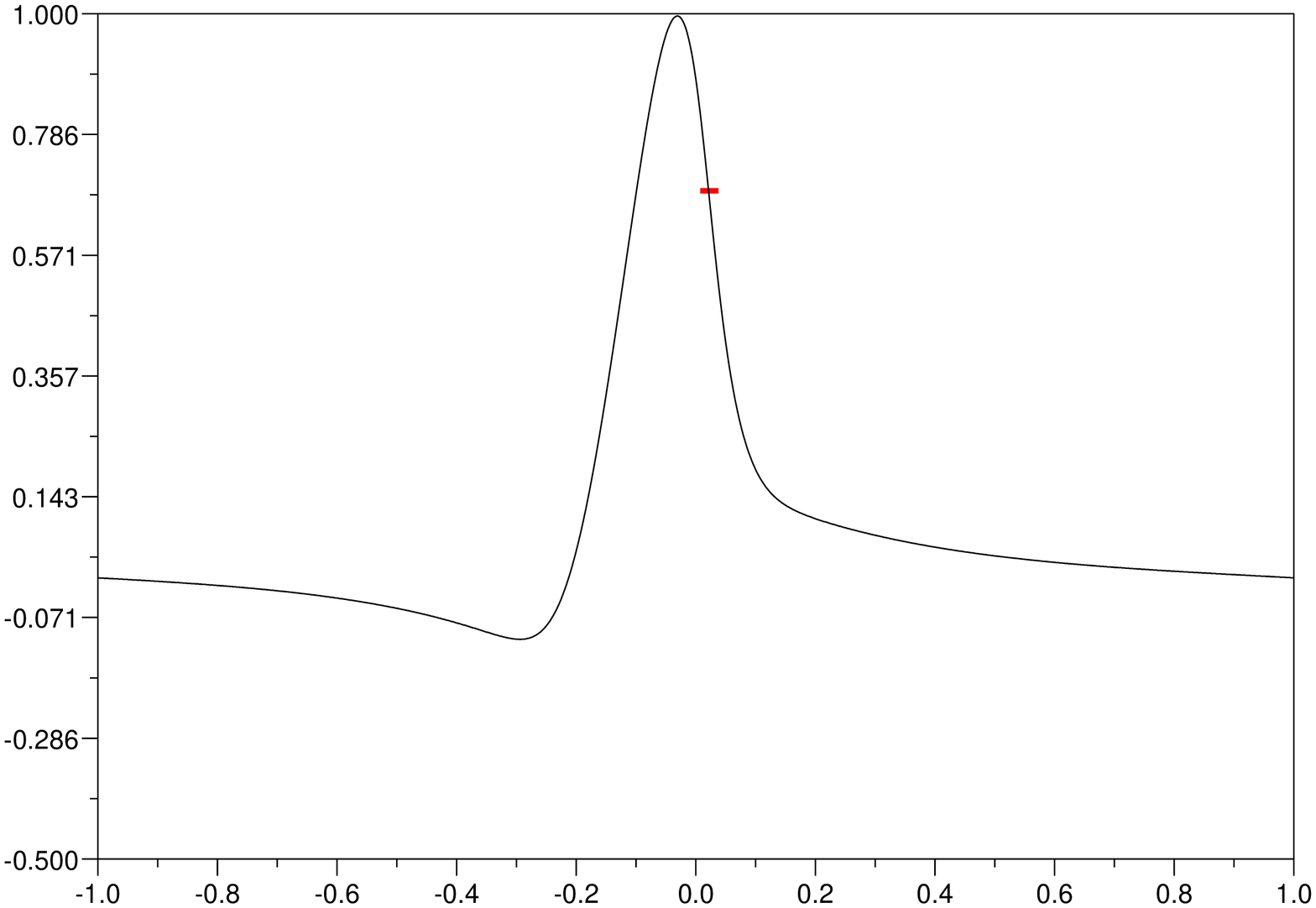}\\
	\includegraphics[width=2.9cm,angle=270]{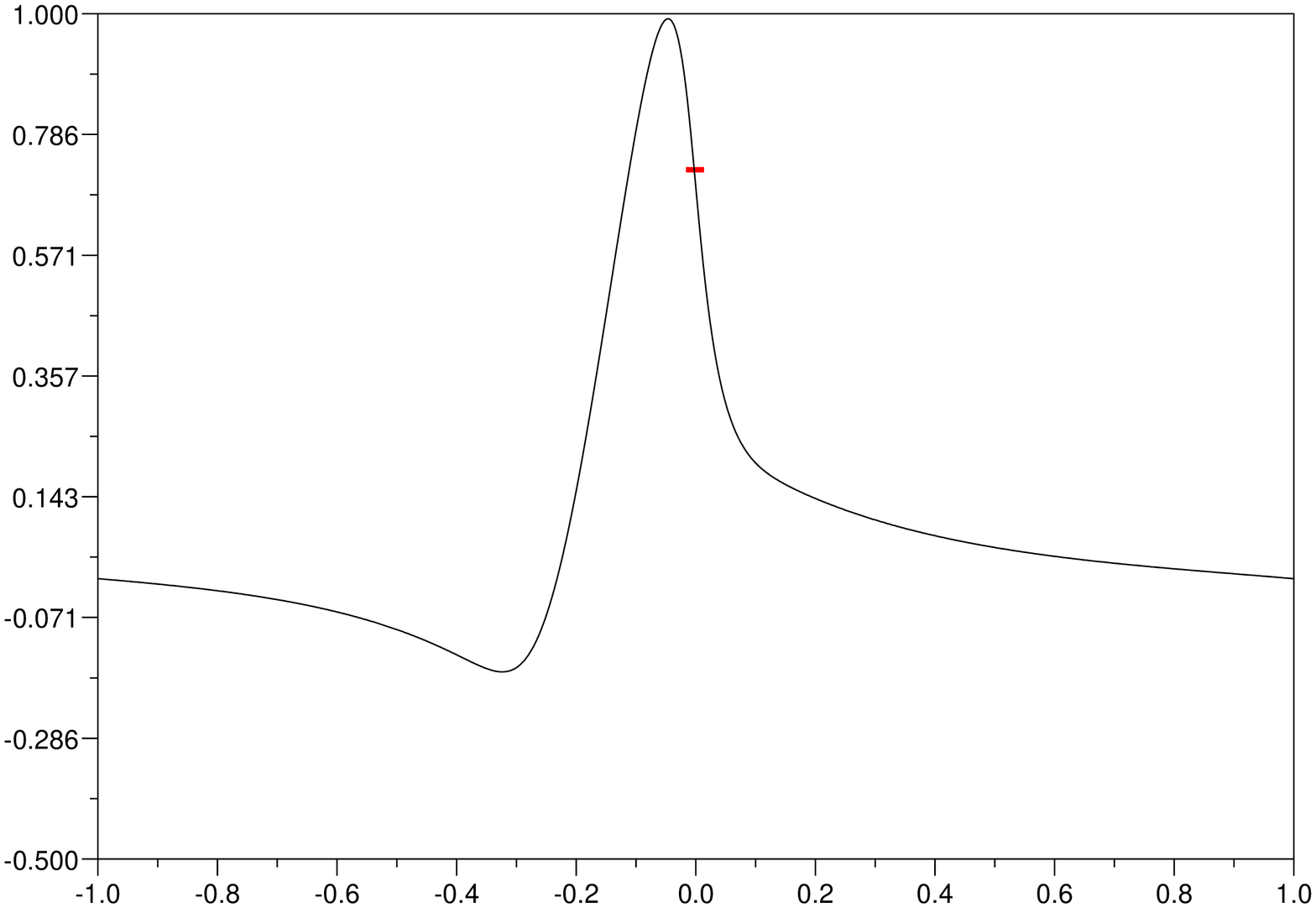}
	\includegraphics[width=2.9cm,angle=270]{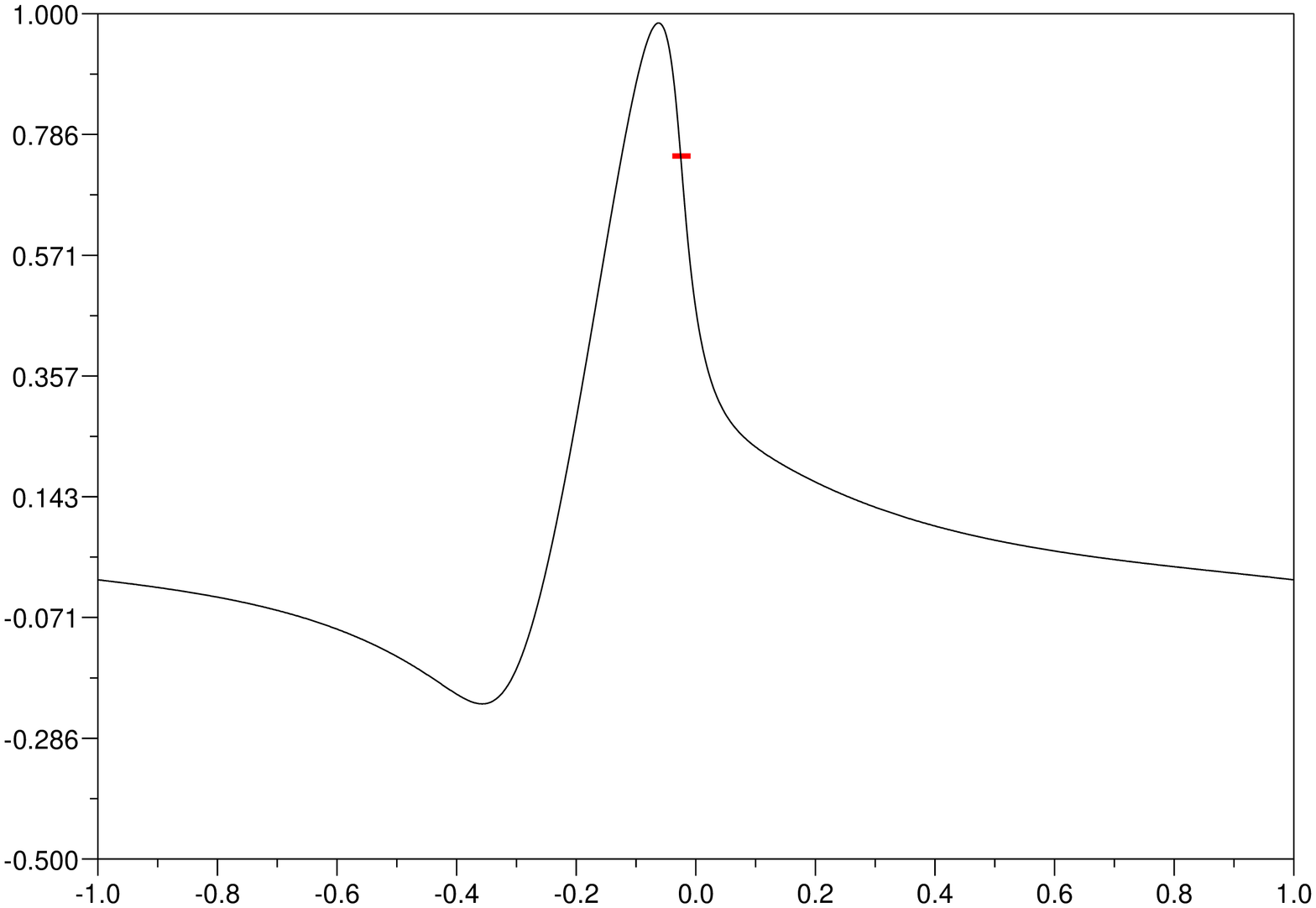}
	\includegraphics[width=2.9cm,angle=270]{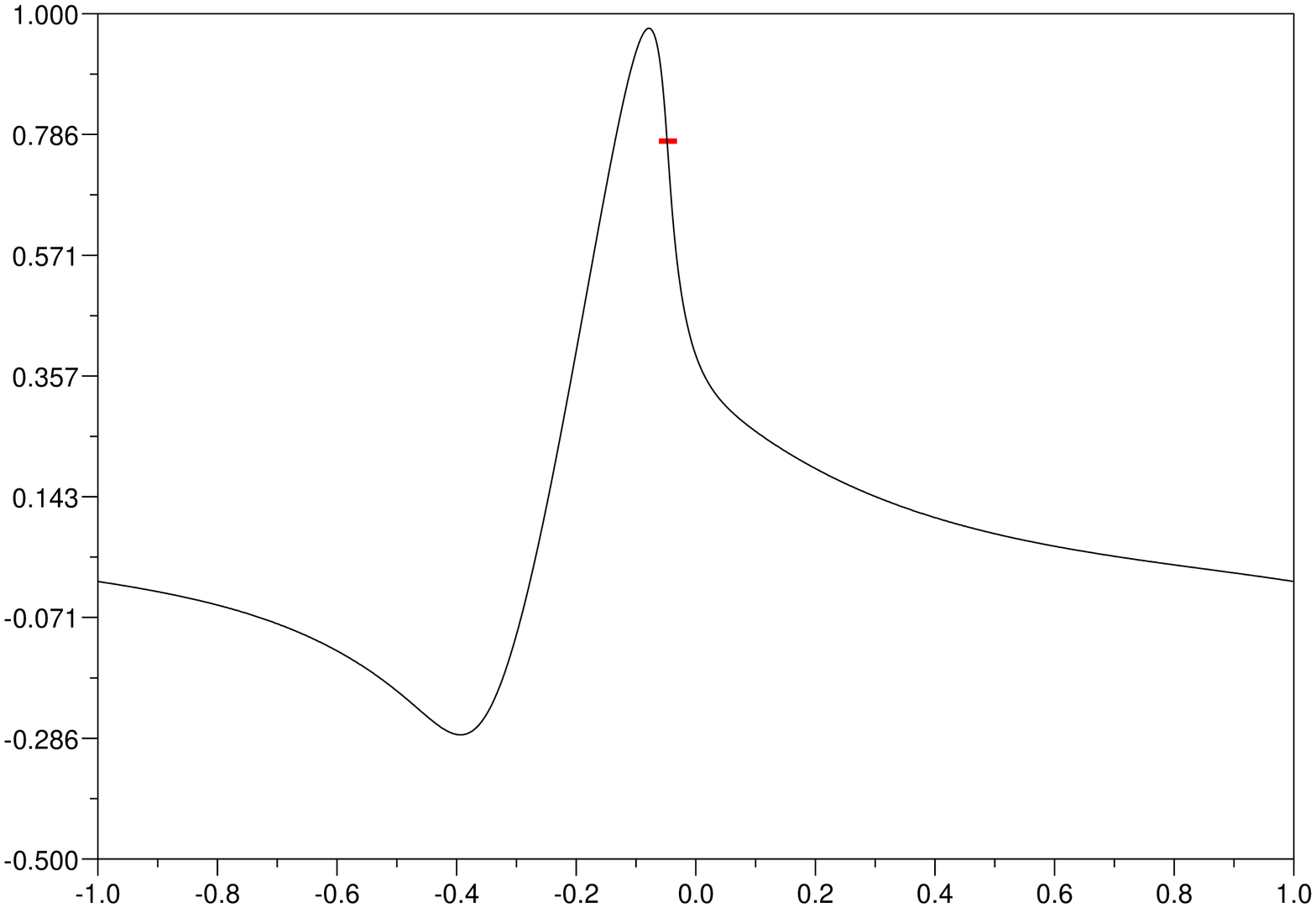}\\
	\includegraphics[width=2.9cm,angle=270]{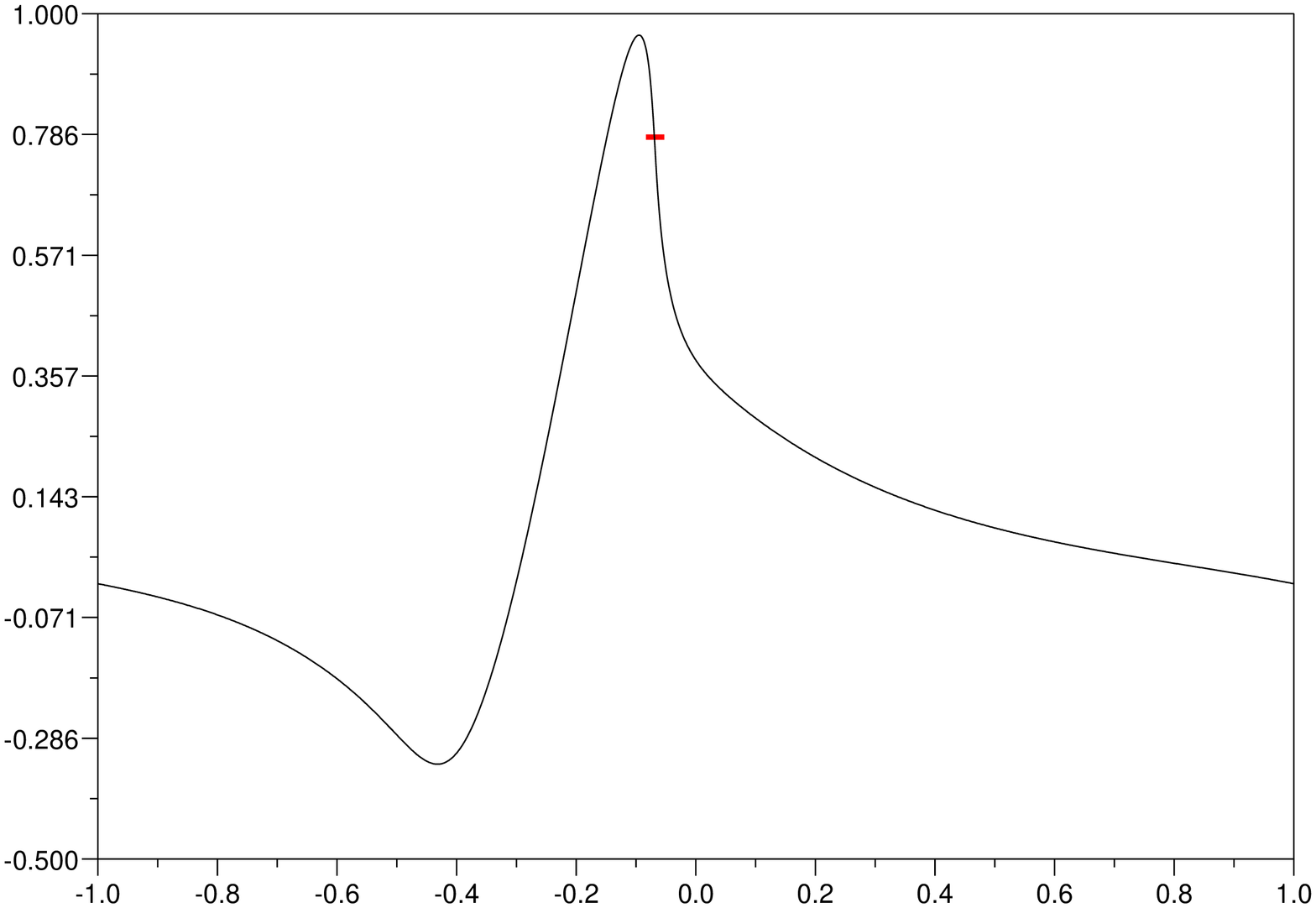}
	\includegraphics[width=2.9cm,angle=270]{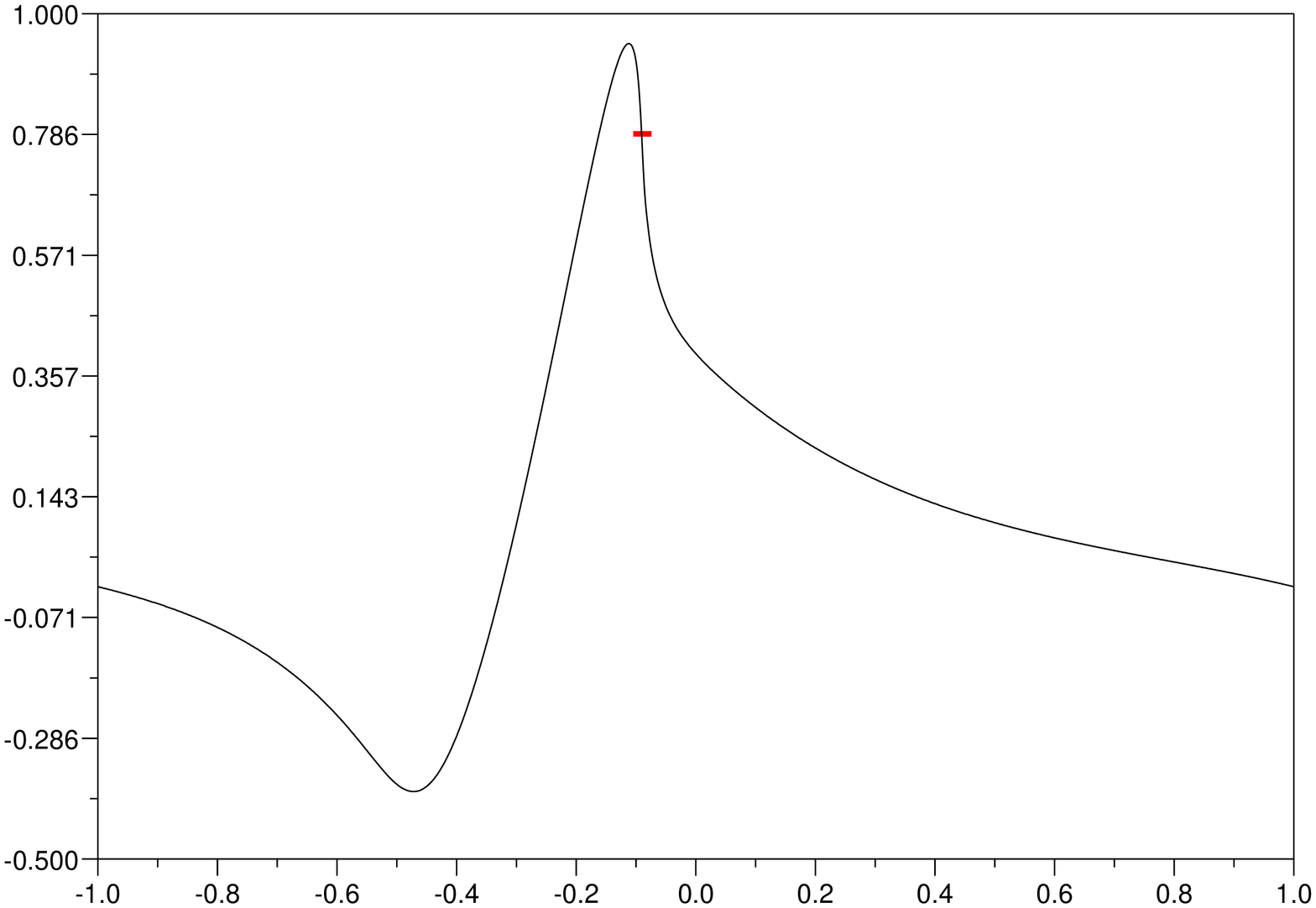}
	\includegraphics[width=2.9cm,angle=270]{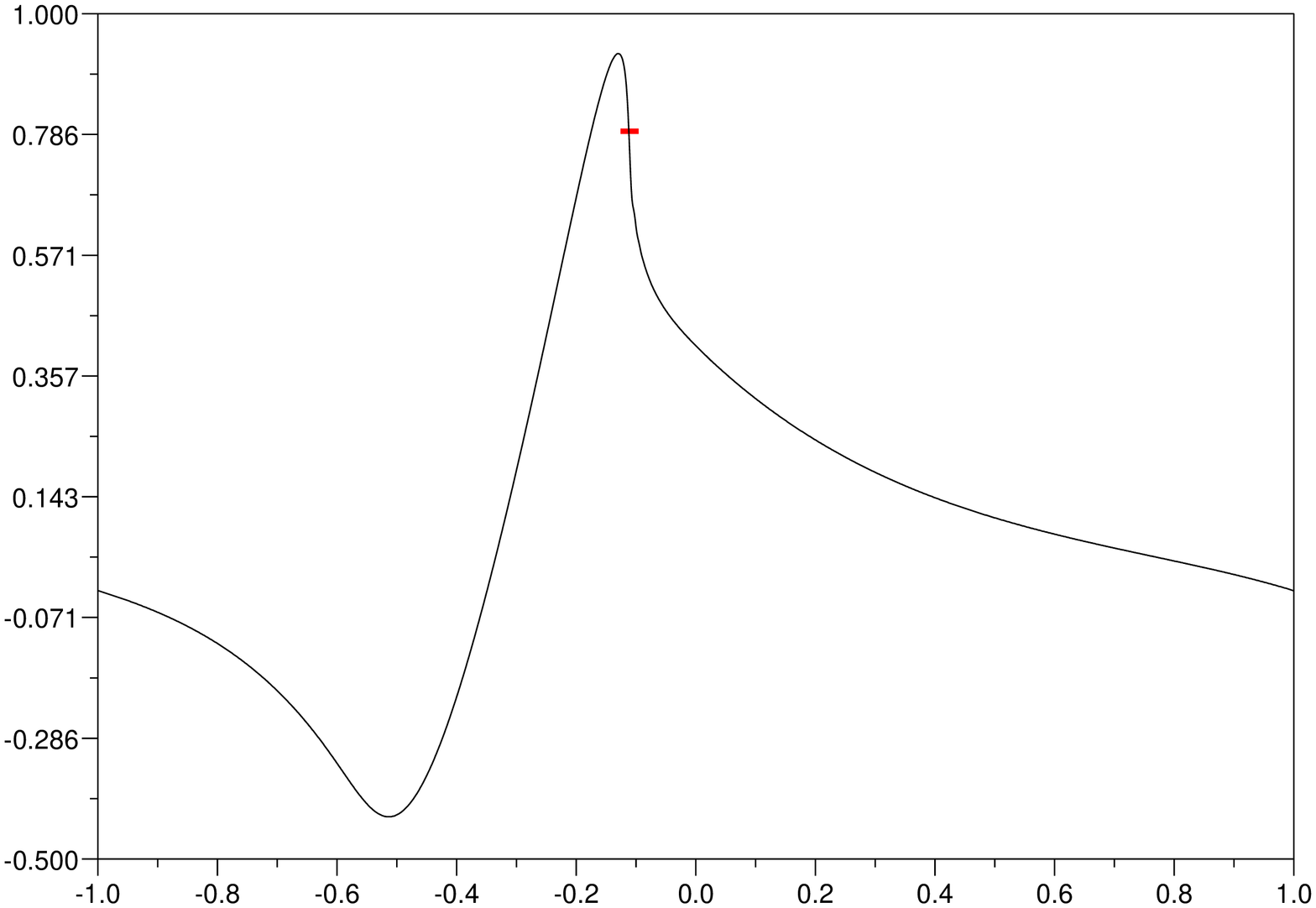}
	\caption{A plunging breaker for the Camassa-Holm equation 
	(\ref{equationCH})}
	\label{figplung}
	\end{center}
\end{figure}

\begin{figure} [htb]
	\begin{center}
	\includegraphics[width=2.9cm,angle=270]{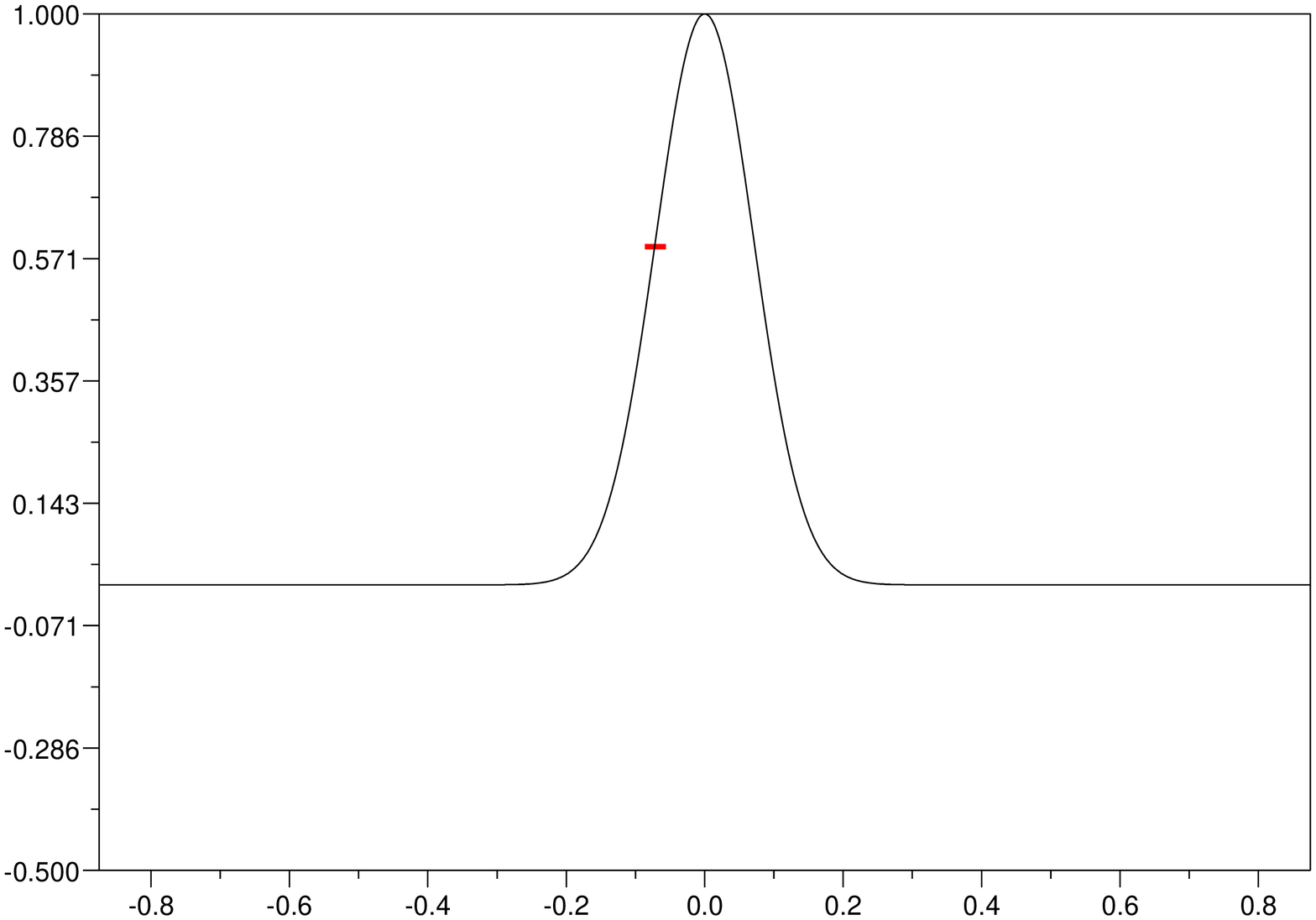}
	\includegraphics[width=2.9cm,angle=270]{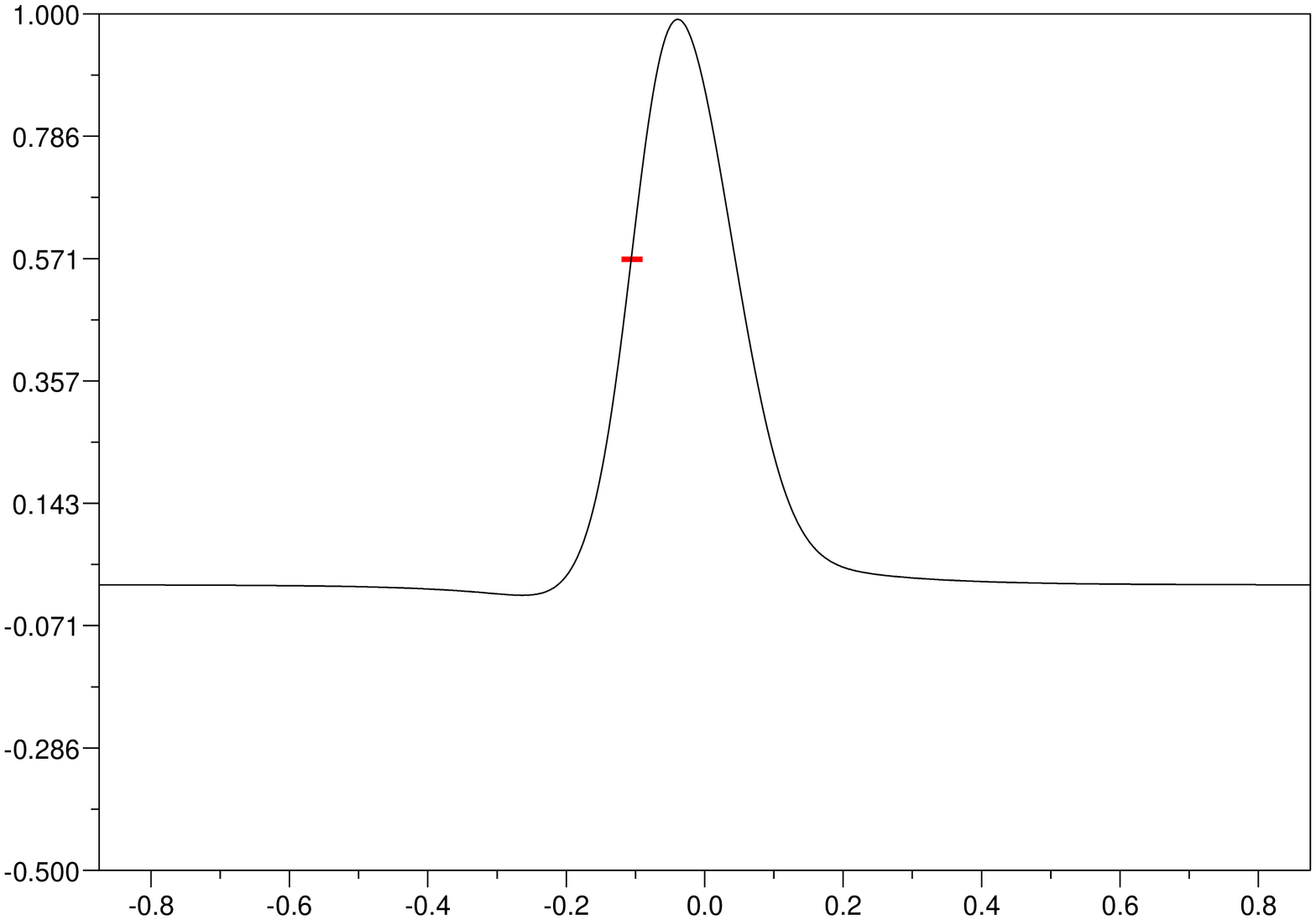}
	\includegraphics[width=2.9cm,angle=270]{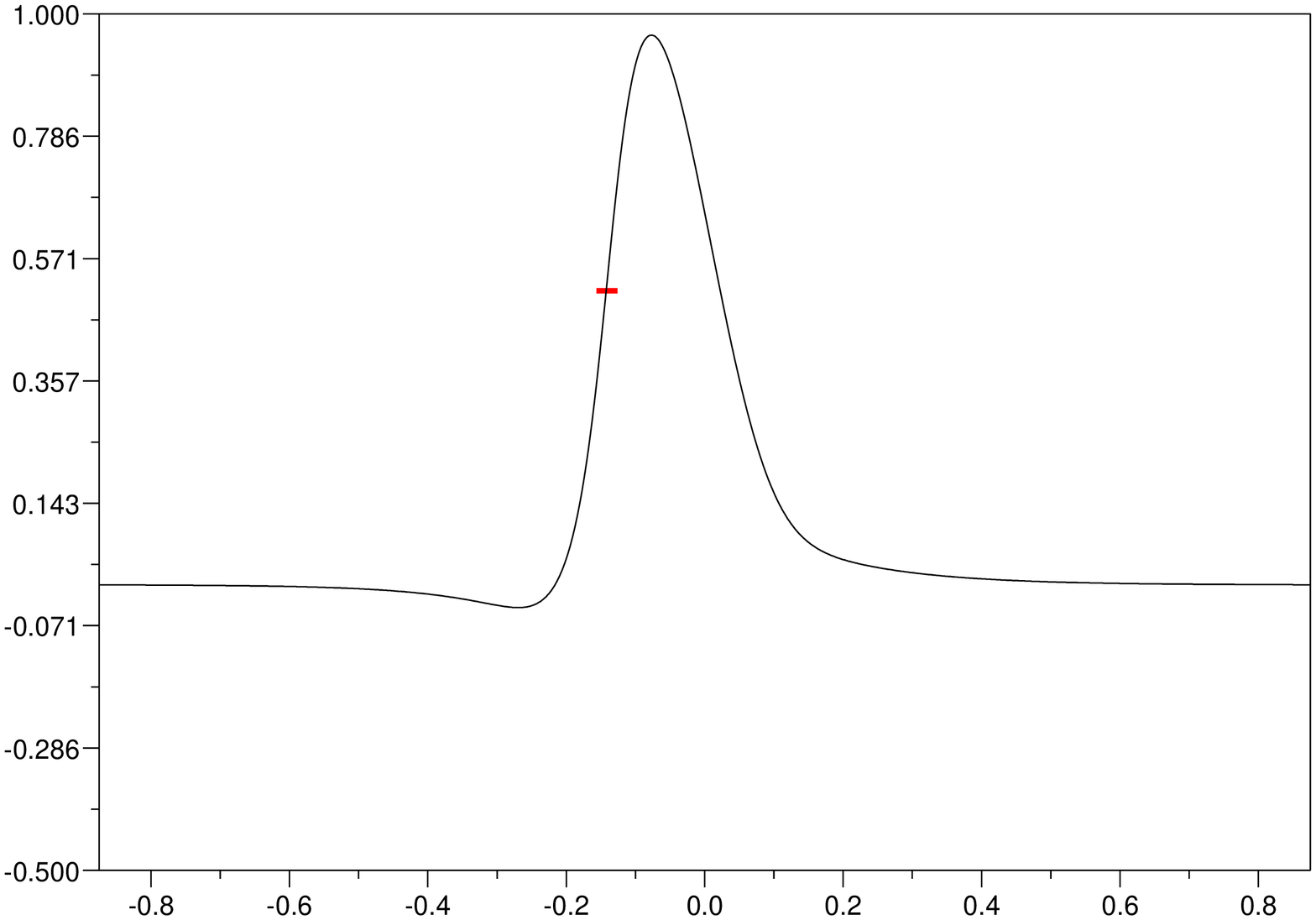}\\
	\includegraphics[width=2.9cm,angle=270]{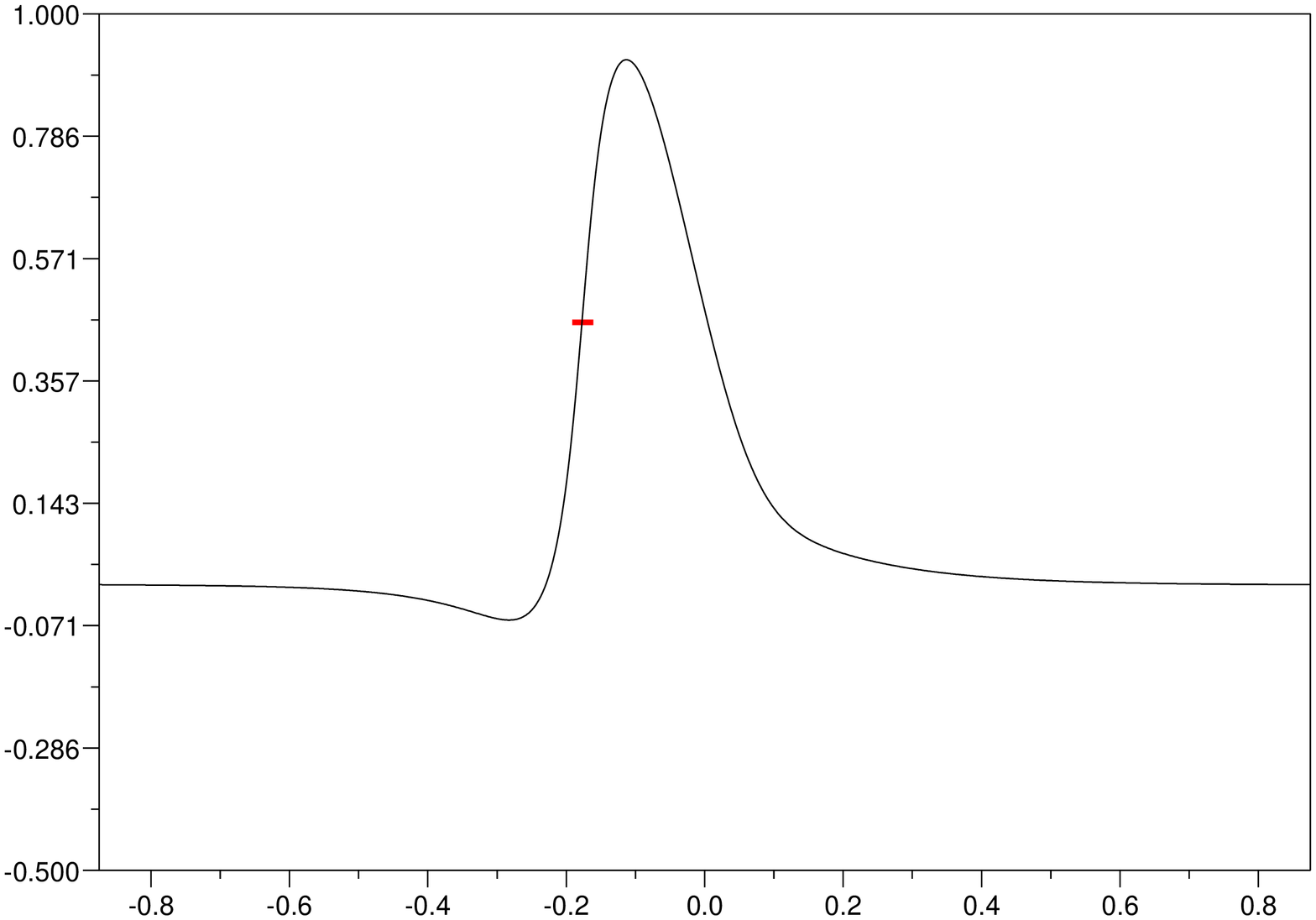}
	\includegraphics[width=2.9cm,angle=270]{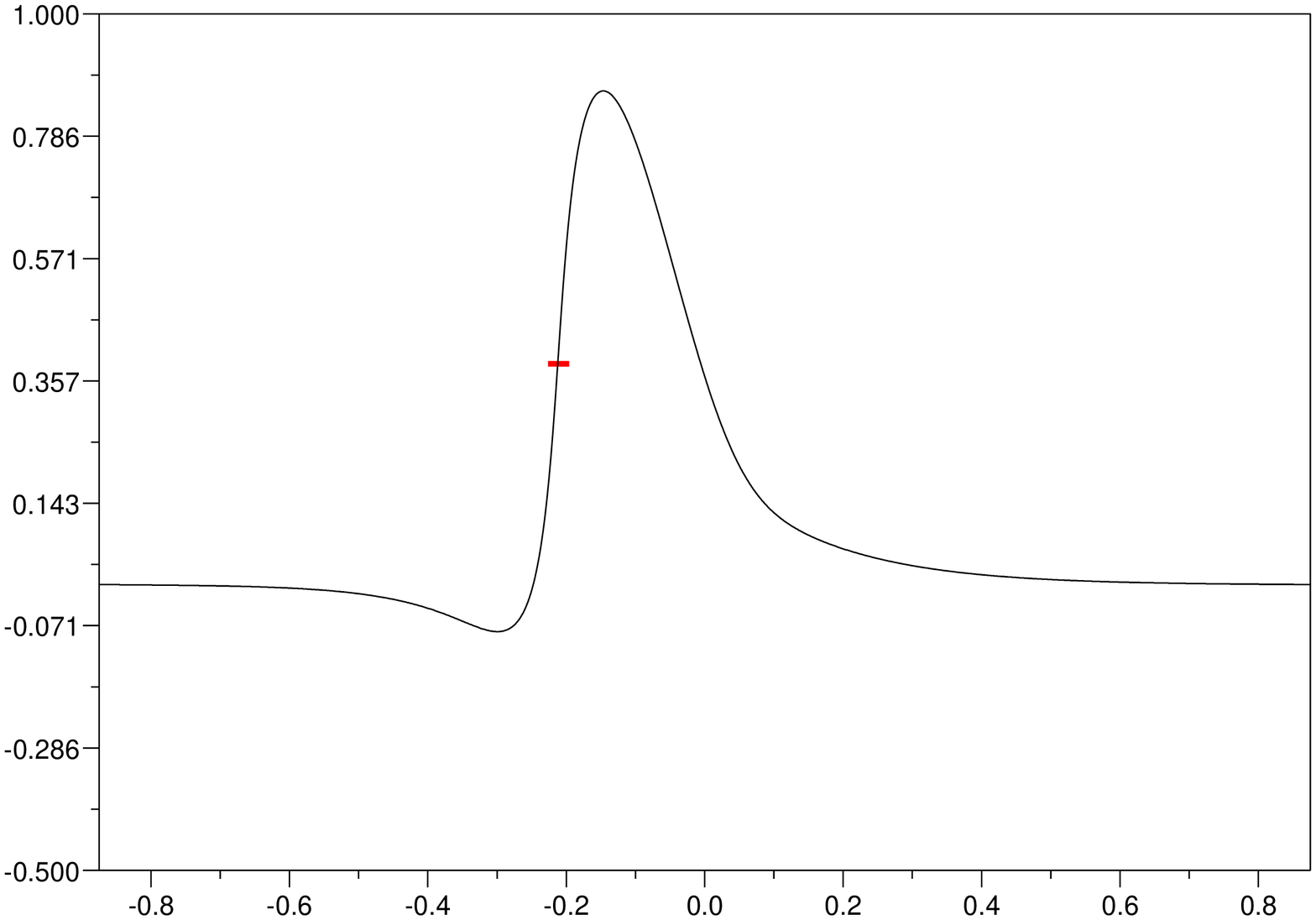}
	\includegraphics[width=2.9cm,angle=270]{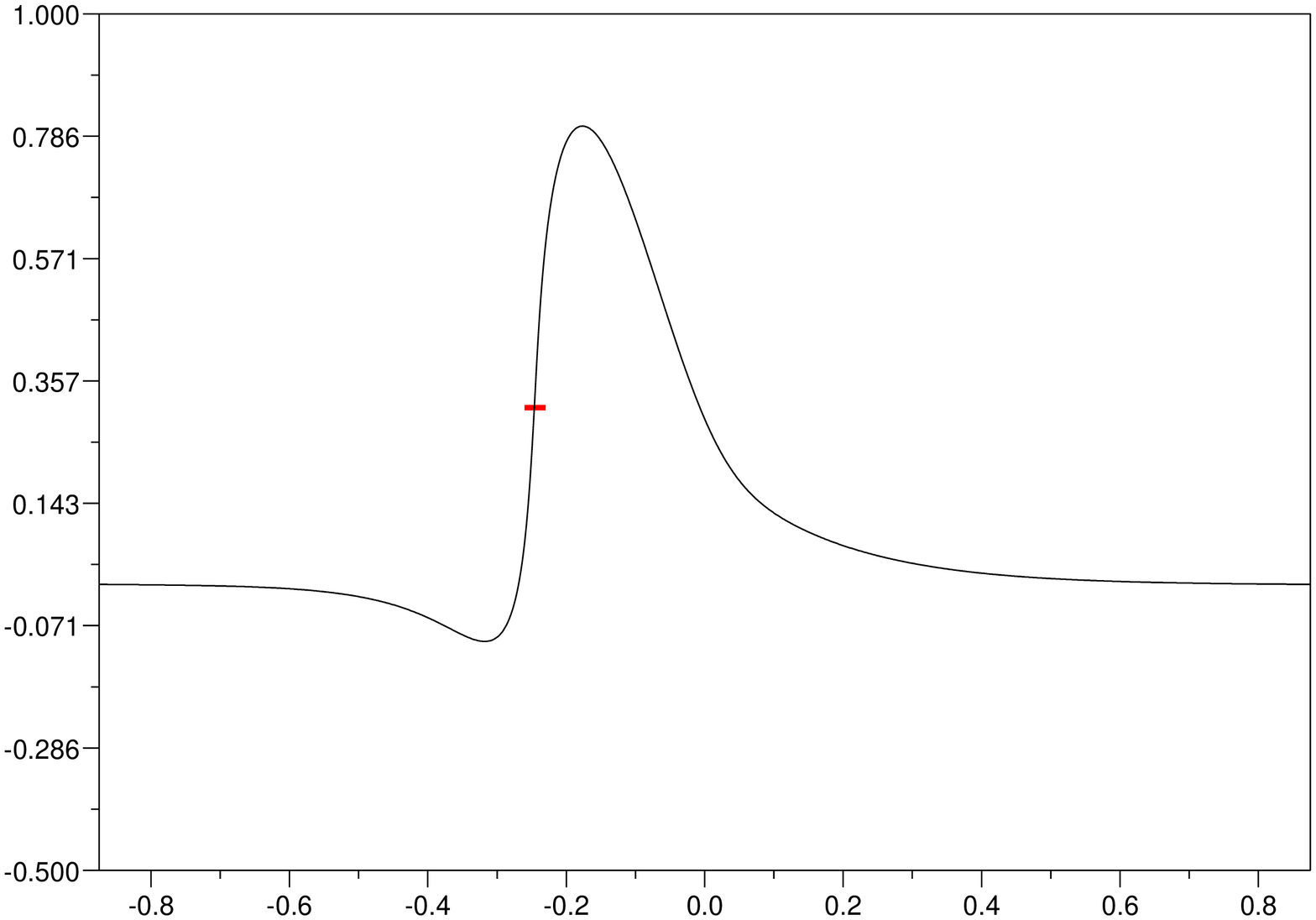}\\
	\includegraphics[width=2.9cm,angle=270]{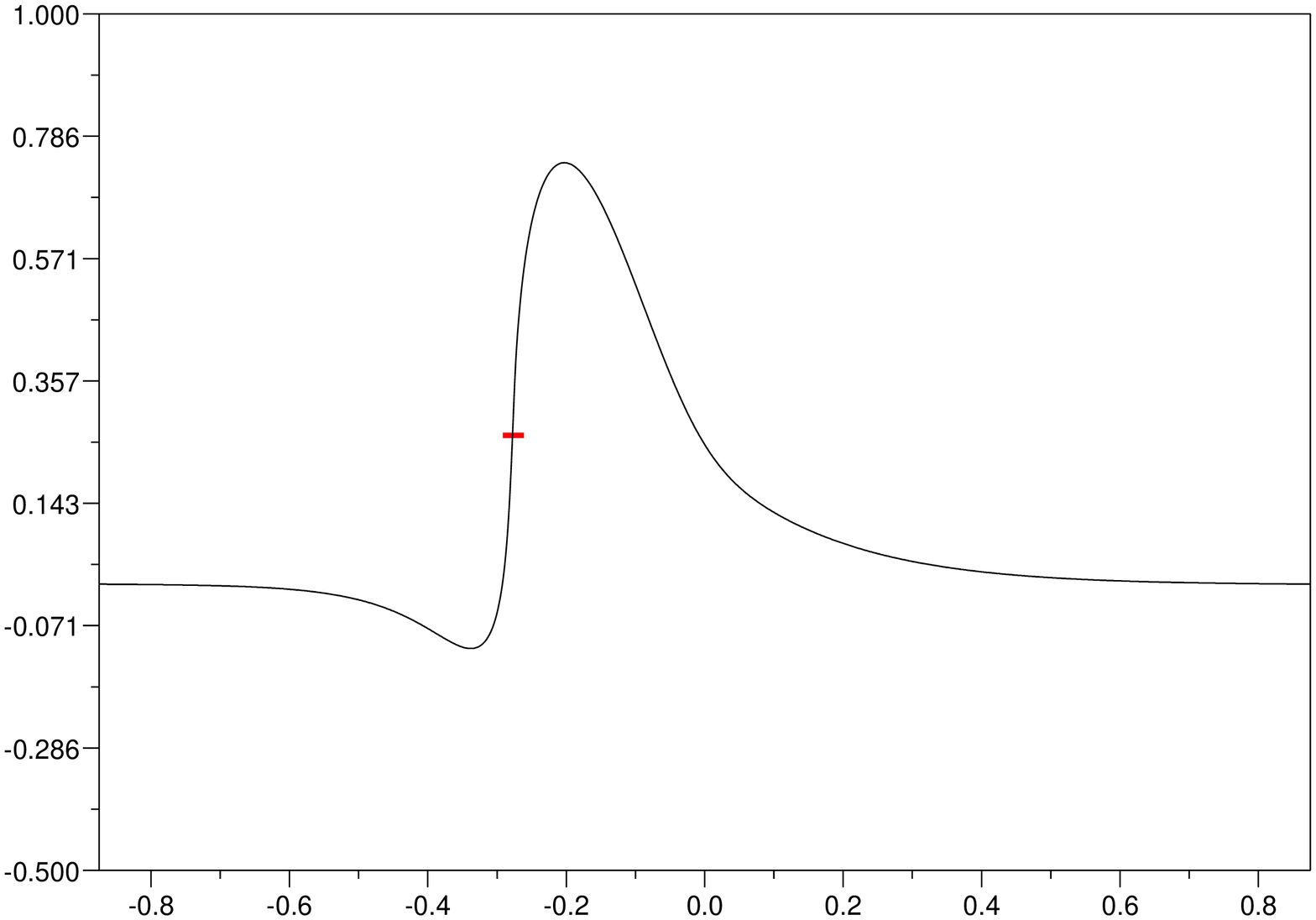}
	\includegraphics[width=2.9cm,angle=270]{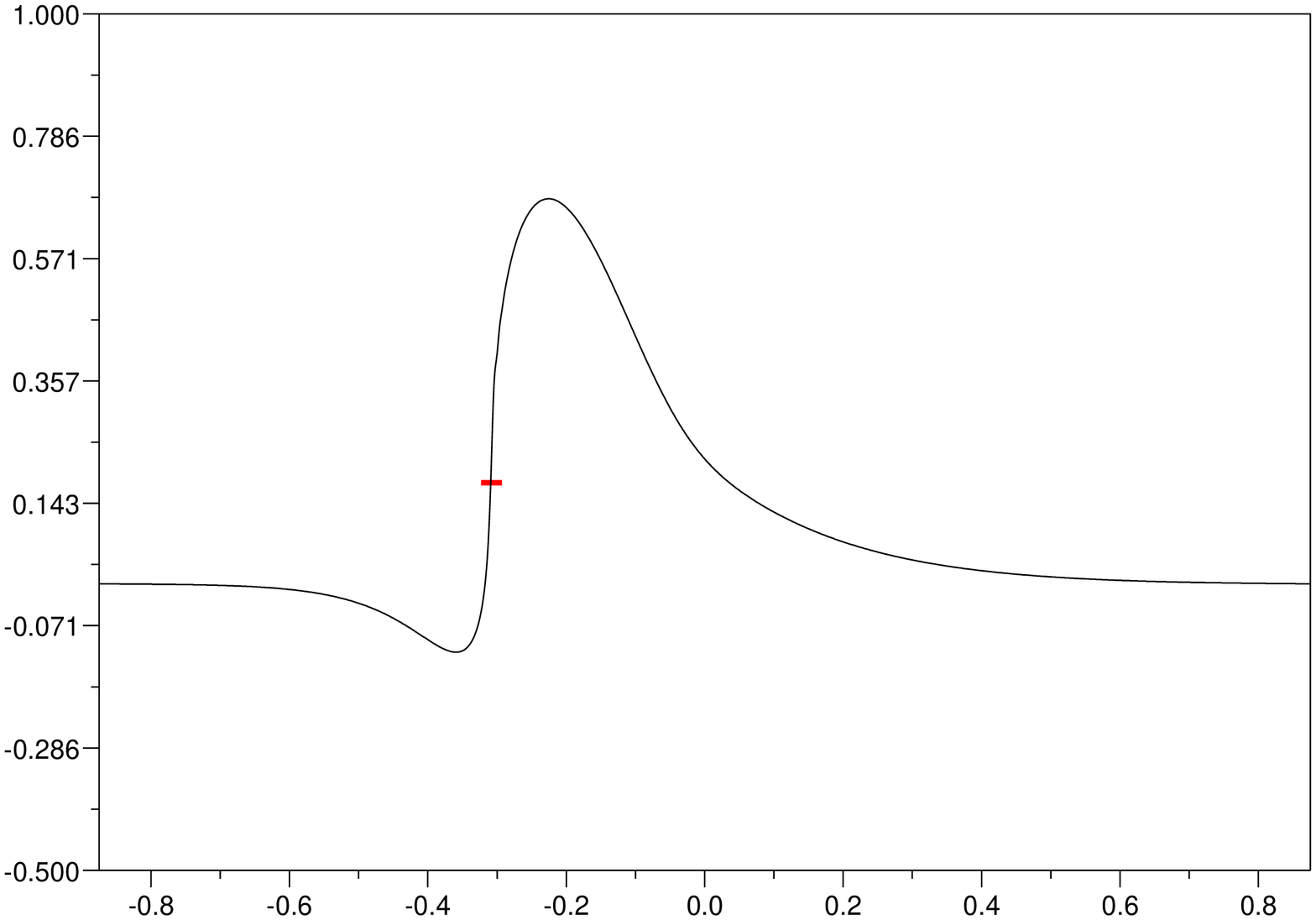}
	\includegraphics[width=2.9cm,angle=270]{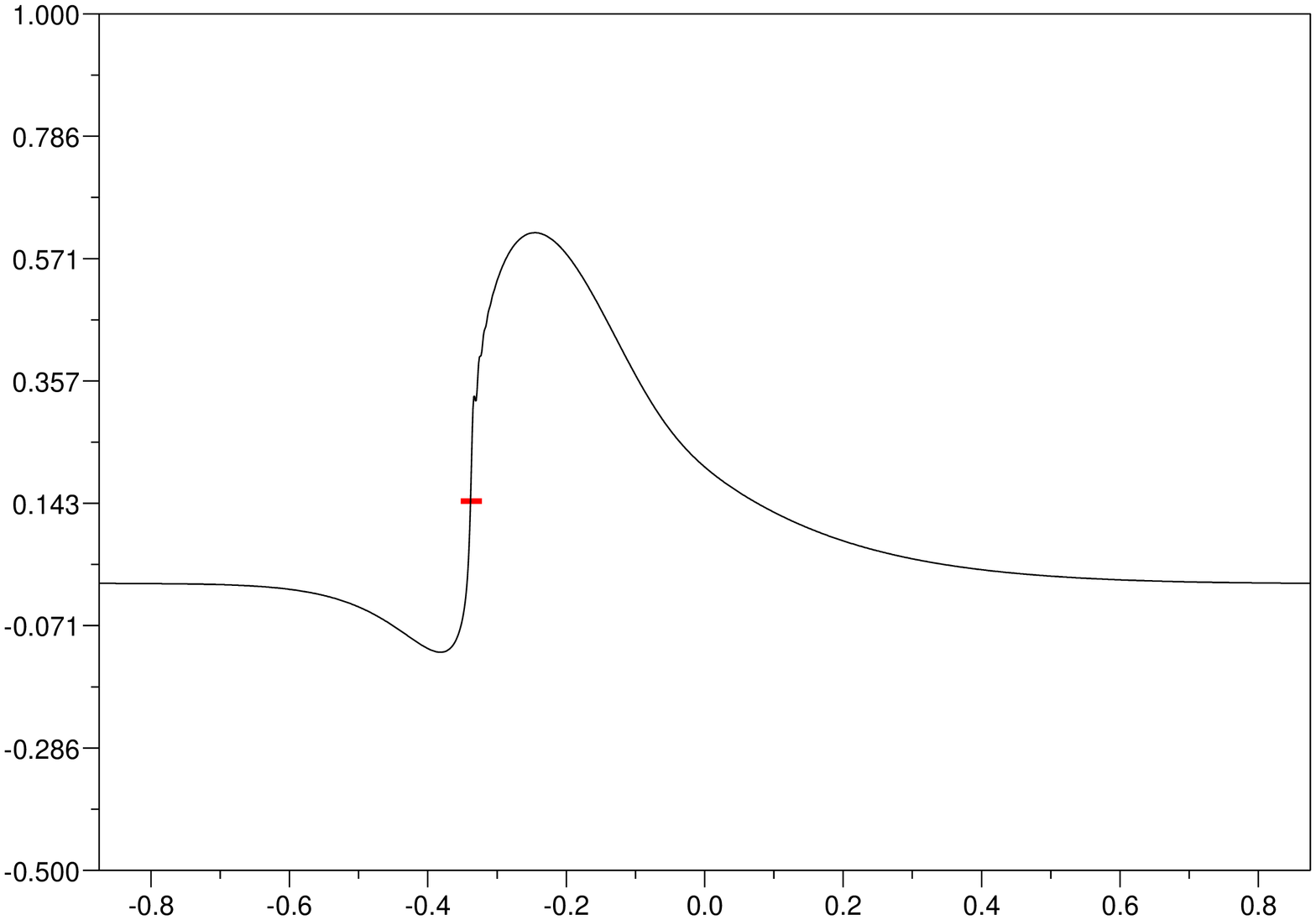}
	\caption{A surging breaker for the surface equation 
	(\ref{eqs3})}
	\label{figsurg}
	\end{center}
\end{figure}

In this section, we use numerical computations to check that the
surface equation (\ref{eqs3}) and the Camassa-Holm equation (\ref{equationCH})
lead respectively to \emph{surging} and \emph{plunging} breakers
as predicted theoretically by Proposition \ref{blowCH2} (and \cite{Co}
for (\ref{equationCH})). We use the same kind of numerical scheme for both
equations; in fact, our scheme works for any equation of the class
(\ref{CHgeneral}) with $\beta<0$.\\
In order to reduce the size of the computational domain, we solve
(\ref{CHgeneral}) in a frame moving at speed $1$, in which
(\ref{CHgeneral}) is replaced by
\begin{eqnarray*}
&u_t+\frac{3}{2}\eps u u_\xi+\eps^2\iota u^2u_\xi+\eps^3\kappa u^3u_\xi +
	\mu ((\alpha-\beta) u_{\xi\xi\xi}+\beta u_{\xi\xi t})\\
&\qquad =\eps\mu(\gamma u u_{\xi\xi\xi}+\delta u_\xi u_{\xi\xi}),\nonumber
\end{eqnarray*}
where $\xi$ stands for $x-t$.\\
The numerical scheme used here is a simple finite difference 
leapfrog/Crank-Nicolson scheme whose semi-discretized version reads
$$
	\forall n\geq 1,\qquad
	(1+\mu\beta\dx^2)\frac{u^{n+1}-u^{n-1}}{2\delta_t}=F[u^n],
$$
with 
$$
	F[u]=-\frac{3}{2}\eps u u_\xi-\eps^2\iota u^2u_\xi-
\eps^3\kappa u^3u_\xi -
	\mu (\alpha-\beta) u_{\xi\xi\xi}
+\eps\mu(\gamma u u_{\xi\xi\xi}+\delta u_\xi u_{\xi\xi}),
$$
 and
where $\delta_t$ is the time step and $u^n\sim u_{\vert_{t=n\delta_t}}$
(to start the induction, that is, for $n=0$, the centered discrete
time derivative must be
replaced by an upwind one).\\
Numerical computations are performed for (\ref{equationCH}) and (\ref{eqs3})
with the same initial value
$$
	u_{\vert_{t=0}}=\exp(-100x^2),
$$
and with $\mu=0.2$, $\eps=\sqrt{\mu}$. Figure \ref{figplung} shows the
formation of a plunging breaker for the solution of (\ref{equationCH}); the
little mark on the curves materializes the point of \emph{minimal} slope. 
For the
same initial data, Figure \ref{figsurg} shows the formation of a surging
breaker for the solution of (\ref{eqs3}); the little mark on the curves
materializes here the point of \emph{maximal} slope.
\providecommand{\href}[2]{#2}

\end{document}